\documentclass[10pt,twoside]{article}
\usepackage[T1]{fontenc}
\usepackage[utf8]{inputenc}
\usepackage{amsmath,amsthm,mathtools}
\usepackage{txfonts}                 
\usepackage{microtype}
\usepackage{enumitem}
\usepackage{titlesec}
\usepackage{fancyhdr}
\usepackage{xcolor}

\usepackage[a4paper,textwidth=140mm,textheight=220mm,headheight=12pt,
            headsep=16pt,footskip=24pt]{geometry}
\definecolor{linkblue}{RGB}{0,70,150}
\usepackage[colorlinks=true,linkcolor=linkblue,citecolor=linkblue,
            urlcolor=linkblue]{hyperref}
\newcommand{\shorttitle}{A logarithmic criterion for first integrals}
\newcommand{\shortauthor}{F.\ Reis, B.\ Sc\'ardua}
\fancypagestyle{firstpage}{\fancyhf{}\fancyfoot[C]{\footnotesize\thepage}}

\titleformat{\section}{\sffamily\bfseries\fontsize{11.5}{14}\selectfont\raggedright}
            {\thesection}{0.6em}{}
\titleformat{\subsection}{\sffamily\bfseries\normalsize\raggedright}
            {\thesubsection}{0.6em}{}
\titlespacing*{\section}{0pt}{24pt plus 4pt minus 2pt}{11pt plus 1pt}
\titlespacing*{\subsection}{0pt}{16pt plus 3pt minus 1pt}{8pt plus 1pt}
\titleformat{\paragraph}[runin]{\normalfont\normalsize\bfseries}{}{0pt}{}
\titlespacing*{\paragraph}{\parindent}{10pt plus 2pt}{0.5em}

\newtheoremstyle{jplain}{9pt plus 2pt minus 1pt}{9pt plus 2pt minus 1pt}
  {\itshape}{0pt}{\sffamily\bfseries}{}{0.6em}
  {\thmname{#1}\thmnumber{ #2}\thmnote{\normalfont\enspace(\textit{#3})}}
\newtheoremstyle{jdefinition}{9pt plus 2pt minus 1pt}{9pt plus 2pt minus 1pt}
  {\normalfont}{0pt}{\sffamily\bfseries}{}{0.6em}
  {\thmname{#1}\thmnumber{ #2}\thmnote{\normalfont\enspace(\textit{#3})}}
\theoremstyle{jplain}
\newtheorem{main}{Theorem}

\newtheorem{maincor}{Corollary}[main]
\newtheorem{theorem}{Theorem}[section]
\newtheorem{lemma}[theorem]{Lemma}
\newtheorem{proposition}[theorem]{Proposition}
\newtheorem{corollary}[theorem]{Corollary}
\theoremstyle{jdefinition}
\newtheorem{definition}[theorem]{Definition}
\newtheorem{example}[theorem]{Example}

\numberwithin{equation}{section}
\makeatletter
\renewcommand{\@biblabel}[1]{#1.}
\renewenvironment{proof}[1][\proofname]{\par\pushQED{\qed}\normalfont
  \topsep6\p@\@plus6\p@\relax\trivlist
  \item[\hskip\labelsep{\sffamily\bfseries\itshape #1}]\hspace{0.35em}\ignorespaces}
  {\popQED\endtrivlist\@endpefalse}
\makeatother

\setlist[enumerate]{leftmargin=2.3em,itemsep=2pt,topsep=4pt}
\setlist[itemize]{leftmargin=1.6em,itemsep=2pt,topsep=4pt}

\newcommand{\C}{\mathbb{C}}
\newcommand{\R}{\mathbb{R}}
\newcommand{\Z}{\mathbb{Z}}
\newcommand{\Q}{\mathbb{Q}}
\newcommand{\N}{\mathbb{N}}
\newcommand{\Pp}{\mathbb{P}}
\newcommand{\Fcal}{\mathcal{F}}
\newcommand{\Gcal}{\mathcal{G}}
\newcommand{\Ocal}{\mathcal{O}}
\newcommand{\tF}{\widetilde{\mathcal{F}}}

\DeclareMathOperator{\ord}{ord}
\DeclareMathOperator{\sing}{sing}
\DeclareMathOperator{\id}{id}
\DeclareMathOperator{\Diff}{Diff}
\DeclareMathOperator{\Hol}{Hol}
\newcommand{\Holv}{\operatorname{Hol}^{\mathrm{virt}}}
\DeclareMathOperator{\hol}{hol}
\DeclareMathOperator{\Inv}{Inv}
\DeclareMathOperator{\Fix}{Fix}
\DeclareMathOperator{\CS}{CS}
\DeclareMathOperator{\Res}{Res}
\DeclareMathAlphabet{\mathcal}{OMS}{cmsy}{m}{n}   
\renewcommand{\Re}{\operatorname{Re}}
\renewcommand{\Im}{\operatorname{Im}}

\hypersetup{
  pdftitle={Growth of leaf counts and a logarithmic criterion for holomorphic first integrals of foliation germs},
  pdfauthor={Fernando Pereira Paulucio Reis, Bruno Sc\'ardua}
}
\begin{document}
\thispagestyle{firstpage}

\begingroup
\raggedright
{\sffamily\bfseries\fontsize{15.5}{19}\selectfont
Growth of leaf counts and a logarithmic criterion for holomorphic first integrals of foliation germs\par}
\vspace{15pt}
{\sffamily\bfseries\fontsize{10.5}{13}\selectfont
 Fernando Reis\textsuperscript{1}\enspace$\cdot$\enspace
 Bruno Sc\'ardua\textsuperscript{2}\par}
\vspace{26pt}
\endgroup
\begingroup
\renewcommand{\thefootnote}{}
\begin{NoHyper}\footnotetext{\hspace*{-1.8em}\begin{tabular}[t]{@{}r@{\ \,}l@{}}
 & Fernando Pereira Paulucio Reis\\
         & \href{mailto:fernandoreis@ipoli.macae.ufrj.br}{fernandoreis@ipoli.macae.ufrj.br}\\[3pt]
   $\boxtimes$       & Bruno Sc\'ardua\\
         & \href{mailto:bruno.scardua@gmail.com}{bruno.scardua@gmail.com}\\[3pt]
 \textsuperscript{1} & Instituto Polit\'ecnico (IPOLI), Universidade Federal do Rio de Janeiro,\\
         & Campus UFRJ--Maca\'e, Av.\ Alu\'izio da Silva Gomes~50, Gl\'oria,\\
         & Maca\'e, RJ 27930-560, Brazil\\[3pt]
 \textsuperscript{2} & Instituto de Matem\'atica, Universidade Federal do Rio de Janeiro (UFRJ),\\
         & CP~68530, Rio de Janeiro, RJ 21945-970, Brazil
 \end{tabular}}\end{NoHyper}
\endgroup

\noindent{\sffamily\bfseries Abstract}\par\nopagebreak\vspace{1pt}
\noindent
For a non-dicritical germ of holomorphic foliation at $0\in\C^2$, we
count the intersections of a leaf with a transverse sector outside radius
$\rho$, taking the supremum over all leaves. A non-constant holomorphic
first integral exists if and only if, at each separatrix, some sector has
sub-logarithmic counts along a sequence $\rho\to0$. For generalized
curves, one separatrix suffices. This gives a sharp gap between bounded
counts, determined by the multiplicities of a primitive first integral,
and counts that grow at least logarithmically. The proof uses a sharp
logarithmic threshold for periodicity in invariance groups of
one-variable equivalence relations. We also give an explicit family of
non-dicritical generalized curves, parametrized by homogeneous quadratic
polynomials, whose leaves are closed off three smooth separatrices with
trivial holonomy; a member of the family admits a non-constant
holomorphic first integral if and only if a single linear condition on
the coefficients holds, the family contains exactly two analytic types,
and in the non-integrable type the leaf counts grow exactly like
$\rho^{-2/3}$. Finally, in the spirit of Suzuki's theorem on
open equivalence relations, a germ admits a non-constant holomorphic
first integral if and only if it is non-dicritical and the closure of the
set of pairs of points on a common leaf is a proper analytic subset.



\vspace{6pt}
\section{Introduction and main results}\label{sec:intro}

Let $\Fcal$ be a germ of singular holomorphic foliation at $0\in\C^2$,
defined by a germ of holomorphic $1$-form $\omega$ with an isolated zero. A
\emph{holomorphic first integral} of $\Fcal$ is a germ of holomorphic
function $H$ with $dH\wedge\omega=0$. Mattei and Moussu characterized the
existence of a non-constant one by a topological condition on the leaves
\cite[Th\'eor\`eme~B]{MM80}. Along a separatrix, a function on a transverse
section extends to a first integral only if it is invariant under the
holonomy of the separatrix \cite[Chapter~II, \S11D]{IY08}. After the
reduction of singularities, however, the holonomy of the separatrices may
miss the returns of the leaves produced near the exceptional divisor. We
measure all these returns by counting the points in which a leaf meets a
sector of a transverse section, and show that the growth of this count
near the vertex decides the existence of a holomorphic first integral,
with a sharp logarithmic threshold. Theorem~\ref{thm:Q52}, like an
earlier example of Moussu \cite{Mou85}, shows that the holonomy of the
separatrices alone cannot decide it.

\paragraph{Sectors and leaf counts.}
A \emph{representative} of $\Fcal$ is the foliation $\Fcal_U$ defined by
$\omega$ on a connected open neighbourhood $U$ of $0$. A \emph{separatrix}
is an irreducible germ $\Gamma$ of invariant analytic curve at $0$, and
$\Fcal_U$ is \emph{adapted} to $\Gamma$ if
$\Gamma^\circ_U=(\Gamma\cap U)\setminus\{0\}$ is a single leaf of
$\Fcal_U$; such representatives exist in every neighbourhood of $0$
(Lemma~\ref{lem:sep}).

\begin{definition}[sectors transverse to a separatrix]\label{def:sector}
Let $\Fcal_U$ be a representative adapted to a separatrix $\Gamma$, and
let $a\in\Gamma^\circ_U$. A \emph{transverse section to $\Gamma$ at $a$}
is a holomorphic disc $\Sigma\subset U$ through $a$, transverse to the
leaves of $\Fcal_U$, with a holomorphic coordinate $z$ on $\Sigma$ such
that $z(a)=0$. For $r>0$ and $\theta_1<\theta_2\leq\theta_1+2\pi$ such that
$\{|z|\leq r\}\subset z(\Sigma)$, the set
\[
S=S(r;\theta_1,\theta_2)=\{x\in\Sigma:\ 0<|z(x)|<r,\ \arg
z(x)\in(\theta_1,\theta_2)+2\pi\Z\}
\]
is a \emph{sector transverse to $\Gamma$}, with vertex $a$ and opening
$\vartheta=\theta_2-\theta_1\in(0,2\pi]$. The coordinate $z$ is part of
the data.
\end{definition}

\begin{definition}[leaf counts]\label{def:counts}
For a sector $S$ transverse to $\Gamma$ in $\Fcal_U$ and $\rho>0$, let
\[
\nu_U(S,\rho)=\sup_L\#\bigl(L\cap S\cap\{|z|>\rho\}\bigr),\qquad
\nu_U(S)=\sup_L\#(L\cap S)=\lim_{\rho\to0}\nu_U(S,\rho),
\]
where $L$ ranges over the leaves of $\Fcal_U$. The sector has
\emph{bounded counts} if $\nu_U(S)<\infty$, and \emph{sub-logarithmic
counts} if $\liminf_{\rho\to0}\nu_U(S,\rho)/\log(1/\rho)=0$. The
\emph{leaf counting function} of $\Gamma$ is
$N_\Gamma(\vartheta)=\inf\nu_U(S)$, $0<\vartheta\leq2\pi$, the infimum
being taken over all representatives adapted to $\Gamma$ and all sectors
of opening $\vartheta$ transverse to $\Gamma$ in them.
\end{definition}

The leaves induce on a transverse section an equivalence relation whose
invariance group is the \emph{virtual holonomy group} of Camacho, Lins Neto
and Sad \cite{CLS92}: the group of germs of biholomorphisms of the section
which preserve the trace of each leaf. Our first result concerns such
groups.

\paragraph{Relations and their invariance groups in dimension one.}
The study of open equivalence relations on analytic spaces goes back to
Cartan \cite{Car60} and Suzuki \cite{Suz74}.
It can be a powerful tool when considering the analytic structure of a holomorphic 
foliation. In this work we pursue this idea, starting with the following.  
For an equivalence relation $\sim$ on a neighbourhood of $0\in\C$, the
\emph{invariance group} $\Inv(\sim)$ is the group of germs
$\gamma\in\Diff(\C,0)$ such that $\gamma(z)\sim z$ for $z$ close to $0$.

\begin{main}\label{thm:onevar}
Let $\sim$ be an equivalence relation on a neighbourhood of $0\in\C$ and
let $S$ be a sector with vertex $0$ and opening $\vartheta\in(0,2\pi]$.
For $\rho>0$, let $\nu(S,\rho)$ be the supremum, over the equivalence
classes $C$, of the number of points of $C\cap S\cap\{|z|>\rho\}$.
\begin{enumerate}[label=\textup{(\alph*)}]
\item If $\nu(S,\rho)\leq N$ for every $\rho>0$, then $\Inv(\sim)$ is a
finite cyclic group of order at most $2\pi N/\vartheta$. The integer
bound $\lfloor2\pi N/\vartheta\rfloor$ is attained for every $N\geq1$ and
every $\vartheta\in(0,2\pi]$.
\item If $\liminf_{\rho\to0}\nu(S,\rho)/\log(1/\rho)=0$, then every
element of $\Inv(\sim)$ has finite order. Consequently, $\Inv(\sim)$ is
isomorphic to a group of roots of unity, and every finitely generated
subgroup of $\Inv(\sim)$ is finite and cyclic.
\item The relation on the unit disc whose classes are $\{0\}$ and the
sets $\{2^kz\}\cap\{|z|<1\}$, $z\neq0$, satisfies
$\nu(S,\rho)\leq\log_2(1/\rho)+1$ for every sector $S$, and $\Inv(\sim)$
contains the germ $z\mapsto z/2$, of infinite order.
\end{enumerate}
\end{main}

Part~(a) is a sectorial and quantitative counterpart of the theorem of
Mattei and Moussu according to which a finitely generated pseudogroup of
$(\C,0)$ with finite orbits is finite \cite[Chapitre~I]{MM80}; see also
\cite{ReRe14}. For non-solvable groups much more is known, since their
orbits are dense in some sectors \cite{She82,Nak94}. Part~(c) shows that
the threshold in (b) is sharp, and the conclusion of (b) cannot be
improved to finiteness: adapting a construction of Rib\'on
\cite[\S7, Remark~25]{Rib26}, we give relations with an infinite
invariance group whose counts grow more slowly than any prescribed
function tending to infinity. Section~\ref{sec:onevar} also determines the
growth produced by a single germ of infinite order according to its
multiplier (Theorem~\ref{thm:trichotomy}), and derives criteria of
finiteness for finitely generated groups of germs
(Corollary~\ref{cor:fg}).

\paragraph{First integrals.}
Another important subject regarding germs of singular holomorphic foliations is the existence of holomorphic first integrals. Our contribution to this subject is  the following criterion.

\begin{main}\label{thm:plane}
Let $\Fcal$ be a non-dicritical germ of singular holomorphic foliation at
$0\in\C^2$.
\begin{enumerate}[label=\textup{(\alph*)}]
\item $\Fcal$ admits a non-constant holomorphic first integral if and only
if every separatrix $\Gamma$ of $\Fcal$ has, in some representative
$\Fcal_U$ adapted to $\Gamma$, a sector $S$ transverse to $\Gamma$ such
that
\[
\liminf_{\rho\to0}\frac{\nu_U(S,\rho)}{\log(1/\rho)}=0 .
\]
\item If $\Fcal$ is a generalized curve, it admits a non-constant
holomorphic first integral as soon as one of its separatrices satisfies
the condition in \textup{(a)}.
\item If $\Fcal$ admits a non-constant holomorphic first integral, then,
for every separatrix $\Gamma$ and every $\vartheta\in(0,2\pi]$,
\[
N_\Gamma(\vartheta)=\Bigl\lceil\frac{m_\Gamma\vartheta}{2\pi}\Bigr\rceil,
\]
where $m_\Gamma$ is the smallest order along $\Gamma$ of a non-constant
holomorphic first integral of $\Fcal$.
\end{enumerate}
\end{main}

In (a) the condition cannot be imposed on one separatrix only, as the weak
separatrix of a saddle-node shows (Example~\ref{ex:necessity}). By (c), a
bound $N$ on a sector of opening $\vartheta$ gives a first integral of
order at most $\lfloor2\pi N/\vartheta\rfloor$ along $\Gamma$, which is
sharp (Example~\ref{ex:sharp}).

\begin{maincor}\label{cor:gap}
Let $\Fcal$ be a non-dicritical generalized curve at $0\in\C^2$ and
$\Gamma$ a separatrix of $\Fcal$. Then exactly one of the following holds.
\begin{enumerate}[label=\textup{(\roman*)}]
\item $\Fcal$ admits a non-constant holomorphic first integral, and
$N_\Gamma(\vartheta)=\lceil m_\Gamma\vartheta/2\pi\rceil<\infty$ for every
$\vartheta\in(0,2\pi]$.
\item Every sector $S$ transverse to $\Gamma$, in every representative
$\Fcal_U$ adapted to $\Gamma$, satisfies
\[
\liminf_{\rho\to0}\frac{\nu_U(S,\rho)}{\log(1/\rho)}>0 .
\]
\end{enumerate}
\end{maincor}

Thus, near a separatrix of a generalized curve, the leaf counts are either
bounded or grow at least logarithmically. Logarithmic growth occurs for
the linear saddles $y\,dx-\alpha x\,dy$, $\alpha\notin\R$
(Example~\ref{ex:log}).

\paragraph{Moderate sectorial first integrals.}
A \emph{moderate sectorial first integral} of $\Fcal_U$ on a sector $S$
transverse to $\Gamma$ is a holomorphic function on $S$, constant on the
trace of each leaf, which admits a non-constant formal power series
$\sum a_jz^j$ as asymptotic expansion on the proper subsectors of $S$
\cite[Definition~1.1]{RS21}. In \cite[Theorem~A]{RS21} such a function is
required on every transverse section, and its expansion is only required
to be nonzero; the second part of Example~\ref{ex:finite} shows that this
statement does not hold as written, and that the expansion must be non-constant. Since the
first non-constant term of the expansion makes the function injective on a
smaller sector (Lemma~\ref{lem:leading}), Theorem~\ref{thm:plane} gives the
corrected statement, with one sector only.

\begin{maincor}\label{cor:modA}
Let $\Fcal$ be a non-dicritical germ of holomorphic foliation at
$0\in\C^2$. If $\Fcal$ is a generalized curve and some separatrix carries
a moderate sectorial first integral, or if every separatrix carries one,
then $\Fcal$ admits a non-constant holomorphic first integral.
\end{maincor}

\paragraph{Closed leaves and periodic holonomy.}
In \cite{Sca16} the second named author studied the germs whose leaves are
\emph{closed off the separatrices}, and asked whether a non-dicritical
generalized curve with this property, each of whose separatrices has
periodic holonomy, necessarily admits a holomorphic first integral
\cite[Question~5.2]{Sca16}. To refute a conjecture of Thom, Moussu
showed that the cuspidal germ $d(y^2+x^3)+x\,(2y\,dx-3x\,dy)$ has a
single separatrix, with trivial holonomy, that all its leaves accumulate
at the origin, and that its projective holonomy group is not abelian
\cite[Proposition~3]{Mou85}; in particular it has no holomorphic first
integral. The closedness of the leaves off the separatrix, which is not
addressed there, follows from the description of the projective holonomy
group in \cite[Exemples~7.16 et~8.24, Proposition~4.21]{Lor21}, so that this
example also answers Question~5.2 in the negative. The following family
has three smooth separatrices and is reduced by one blow-up; the
closedness of the leaves, the existence of a first integral and the
analytic type are all decided explicitly.

\begin{main}\label{thm:Q52}
Let $F=xy(y-x)$, let $\Delta=a_0x^2+a_1xy+a_2y^2$ be a homogeneous
quadratic polynomial, and let $\Fcal_\Delta$ be the germ at $0\in\C^2$ of
the foliation defined by
\[
\omega_\Delta=\tfrac13\,dF+\Delta\,(x\,dy-y\,dx).
\]
Let $U$ be a connected open neighbourhood of the origin of\/ $\C^2$ in
which $\omega_\Delta$ has no zero other than $0$, and let $\Fcal_U$ be
the corresponding representative of $\Fcal_\Delta$.
\begin{enumerate}[label=\textup{(\alph*)}]
\item $\Fcal_\Delta$ is a non-dicritical generalized curve; one blow-up
reduces it to three analytically linearizable resonant saddles, with
eigenvalue ratio $-3$, its separatrices are the three lines of
$\{F=0\}$, and each of them has trivial holonomy.
\item Every leaf of $\Fcal_U$ not contained in $\{F=0\}$ is a closed
subset of $U\setminus\{F=0\}$: the leaves of $\Fcal_\Delta$ are closed
off the separatrices.
\item $\Fcal_\Delta$ admits a non-constant holomorphic first integral if
and only if $2a_0+a_1+2a_2=0$. In this case
$3\Delta=a(x^2-2xy)+b(y^2-2xy)$ for unique $a,b\in\C$, and
$F/(1-ax-by)^3$ is a first integral.
\item $\Fcal_\Delta$ is analytically conjugate to the germ defined by $dF$
if $2a_0+a_1+2a_2=0$, and to the germ defined by
$\frac13\,dF+xy\,(x\,dy-y\,dx)$ otherwise. More precisely, if
$\Delta'=\alpha\Delta+\frac13\bigl(p(x^2-2xy)+q(y^2-2xy)\bigr)$ with
$\alpha\in\C^*$ and $p,q\in\C$, then the germ of biholomorphism
$(x,y)\mapsto(x,y)/(\alpha+px+qy)$ conjugates $\Fcal_\Delta$ to
$\Fcal_{\Delta'}$.
\item Let $2a_0+a_1+2a_2\neq0$, let $\Gamma$ be a separatrix of
$\Fcal_\Delta$, $\Fcal_V$ a representative adapted to $\Gamma$, and $S$ a
sector transverse to $\Gamma$ in it. There are $c>0$ and $\rho_0>0$ such
that $\nu_V(S,\rho)\geq c\,\rho^{-2/3}$ for $0<\rho<\rho_0$; if the
radius of $S$ is small enough, there is also $C>0$ such that
$\nu_V(S,\rho)\leq C\rho^{-2/3}$ for $0<\rho<\rho_0$.
\end{enumerate}
\end{main}

Thus the family contains exactly two analytic types, and its parameters
are not moduli. In Moussu's cuspidal setting, by contrast, \'Ecalle and
Yoccoz found infinitely many analytic types whose separatrix has trivial
holonomy \cite[p.~164]{Mou85}.

When $2a_0+a_1+2a_2\neq0$, the obstruction lies in the projective holonomy
group. We show in Section~\ref{sec:Q52} (Theorem~\ref{thm:group}) that it
is analytically conjugate to the group generated by $\zeta^2x$ and
$x/(\zeta+\lambda x)$, $\zeta=e^{2\pi i/3}$, with
$\lambda=2a_0+a_1+2a_2$: a crystallographic group
$\Z^2\rtimes\Z/3\Z$, with orbits closed off the origin, if $\lambda\neq0$,
and a cyclic group of order three if $\lambda=0$. In Moussu's example the
corresponding group has a rotation of order six \cite[Exemple~8.24]{Lor21}.
By part~(e), for $\lambda\neq0$ the leaf counts grow at least like
$\rho^{-2/3}$ on every sector transverse to each of the three
separatrices, and exactly like $\rho^{-2/3}$ on every sufficiently small
one: the exponent is the rank $2$ of the lattice of translations
of the projective holonomy group, divided by the multiplicity $3$ of the
first integral $\tau\xi^3$ at the corners (Section~\ref{sec:growthC}).
This is far above the logarithmic threshold of
Corollary~\ref{cor:gap}, although the holonomy of these separatrices is
trivial: the growth of the counts along a separatrix is produced by the
holonomy of the exceptional divisor, which the holonomy of the separatrix
does not see.

\paragraph{Analytic graphs.}
Suzuki proved that an open equivalence relation on a reduced complex space
has an analytic quotient if and only if its graph is analytic
\cite[Th\'eor\`eme~1]{Suz74}. The following notion is the local analogue
for germs of foliations.

\begin{definition}[graph of a foliation]\label{def:graph}
The \emph{graph} of a representative $\Fcal_U$ is the closure
$\mathcal G(\Fcal_U)$, in $U\times U$, of the set of pairs of points of
$U\setminus\sing(\Fcal)$ which lie on a common leaf of $\Fcal_U$. The germ
$\Fcal$ has an \emph{analytic graph} if, for some representative, this
set is an analytic subset of $U\times U$ different from $U\times U$.
\end{definition}

\begin{main}\label{thm:graph}
A germ $\Fcal$ of singular holomorphic foliation at $0\in\C^2$ admits a
non-constant holomorphic first integral if and only if it is
non-dicritical and has an analytic graph. In this case $\Fcal$ is a
generalized curve, and, for every primitive first integral $F_0$ of
$\Fcal$, there are representatives $\Fcal_V$ in every neighbourhood of
the origin such that
$\mathcal G(\Fcal_V)=\{(p,q)\in V\times V:\ F_0(p)=F_0(q)\}$.
\end{main}

If $F_0$ is a primitive first integral, the graph is the analytic set
$\{F_0(p)=F_0(q)\}$ for suitable representatives. Conversely, an analytic
graph forces bounded counts on small sectors transverse to every
separatrix, so that the reverse implication follows from
Theorem~\ref{thm:plane}(a). In Suzuki's terms, the relation induced near
the base point of a transverse section to $\Gamma$ has the analytic
quotient given by $F_0$, of local degree $m_\Gamma$. When
$2a_0+a_1+2a_2\neq0$, the germ $\Fcal_\Delta$ of Theorem~\ref{thm:Q52} has
leaves closed off the separatrices but a graph which is not analytic.

\paragraph{Structure of the paper.}
This paper is organized as follows. Section~\ref{sec:prelim} collects the preliminaries. Section~\ref{sec:onevar}
contains the proof of Theorem~\ref{thm:onevar} and the results on orbit
growth of germs and groups. Section~\ref{sec:foliations} contains the
proofs of Theorem~\ref{thm:plane}, of its corollaries and the examples.
Section~\ref{sec:Q52} contains the proof of Theorem~\ref{thm:Q52},
including the growth of its leaf counts, and the study of the projective holonomy group (Theorem~\ref{thm:group}),
and Section~\ref{sec:graph} that of
Theorem~\ref{thm:graph}.
\section{Preliminaries}\label{sec:prelim}

\subsection{Leaves, transverse sections and holonomy}\label{sec:basics}

A germ of singular holomorphic foliation $\Fcal$ at $0\in\C^2$ is defined
by a germ of holomorphic $1$-form $\omega$ with an isolated zero at $0$,
unique up to a unit. A \emph{representative} $\Fcal_U$ is defined by
$\omega$ on a connected open neighbourhood $U$ of $0$; its \emph{leaves}
are the maximal connected integral curves of $\ker\omega$ in
$U\setminus\sing(\Fcal)$, and we write $x\sim_Uy$ when $x$ and $y$ lie on
the same leaf. If $U'\subset U$, every leaf of $\Fcal_{U'}$ is contained in
a leaf of $\Fcal_U$. Every regular point has a \emph{flow box}: a connected
neighbourhood $V$ and a submersion $f\colon V\to\C$ with connected fibres,
the \emph{plaques}, such that $\omega=\upsilon\,df$ with $\upsilon$ a unit.
A \emph{transverse section} is an embedded holomorphic disc on which
$\omega$ vanishes nowhere; at a point $a$ we fix on it a coordinate $z$
with $z(a)=0$, which identifies the germs of biholomorphisms of
$(\Sigma,a)$ with $\Diff(\C,0)$. Following the plaques along a path $c$ in
a leaf from $a$ to $b$ gives the \emph{holonomy map}
$\hol_c\colon(\Sigma_a,a)\to(\Sigma_b,b)$, which depends only on the
homotopy class of $c$ and satisfies $\hol_c(x)\sim_Ux$; the holonomy maps
of the loops at $a$ form the holonomy group of the leaf
\cite[Chapter~IV]{CLN85}.

\subsection{Virtual holonomy and invariance groups}\label{sec:virtual}

By construction, the maps in the holonomy group preserve the trace of each
leaf on the transverse section. This property is shared by a larger group,
introduced in \cite[Definition~2, p.~432]{CLS92}; see also
\cite{CaSc00} and \cite[\S3.2]{RS21}.

\begin{definition}\label{def:G}
Let $\sim$ be an equivalence relation on an open neighbourhood of $0$ in
$\C$. The \emph{invariance group} of $\sim$ is the set $\Inv(\sim)$ of the
germs $\gamma\in\Diff(\C,0)$ such that $\gamma(z)\sim z$ for every $z$ in
some neighbourhood of $0$. Let $\Fcal_U$ be a representative of $\Fcal$ and
$\Sigma$ a transverse section at $a$. The \emph{virtual holonomy group} of
$\Fcal_U$ with respect to $\Sigma$ and $a$ is
\[
\Holv(\Fcal_U,\Sigma,a)=\{\gamma\in\Diff(\Sigma,a):\ L_{\gamma(x)}=L_x
\text{ for all }x\in(\Sigma,a)\},
\]
where $L_x$ denotes the leaf of $\Fcal_U$ through $x$. In other words,
$\Holv(\Fcal_U,\Sigma,a)=\Inv(\sim)$, where $\sim$ is the restriction of
$\sim_U$ to $\Sigma$, read in the coordinate $z$.
\end{definition}

For the relation defined by the level sets of a function, $\Inv(\sim)$ is
the invariance group of that function in the sense of
\cite[Chapitre~I, \S1]{MM80} and \cite[\S2.4]{RS21}. The same definitions
apply to a holomorphic foliation $\Gcal$ with isolated singularities on a
complex surface $X$, with respect to the leaves of $\Gcal$ in
$X\setminus\sing(\Gcal)$; we then write $\Holv(\Gcal,\Sigma,a)$. The
following properties are immediate consequences of the definitions.

\begin{lemma}\label{lem:vir}
\begin{enumerate}[label=\textup{(\alph*)}]
\item $\Inv(\sim)$ is a subgroup of $\Diff(\C,0)$.
\item If $\sim'$ is an equivalence relation on a neighbourhood of $0$ such
that $z\sim'w$ implies $z\sim w$, then $\Inv(\sim')\subset\Inv(\sim)$. In
particular, $\Holv(\Fcal_{U'},\Sigma,a)\subset\Holv(\Fcal_U,\Sigma,a)$ if
$\Sigma\subset U'\subset U$.
\item $\Hol(L,\Sigma,a)\subset\Holv(\Fcal_U,\Sigma,a)$, where $L$ is the
leaf of $\Fcal_U$ through $a$.
\item Let $\Sigma,\Sigma'$ be transverse sections at $a,a'$, and let
$\sigma\colon(\Sigma,a)\to(\Sigma',a')$ be a germ of biholomorphism such
that $\sigma(x)\sim_U x$ for $x$ close to $a$; for instance, the holonomy
map of a path in a leaf of $\Fcal_U$. Then
$\Holv(\Fcal_U,\Sigma',a')=\sigma\circ\Holv(\Fcal_U,\Sigma,a)\circ\sigma^{-1}$.
\end{enumerate}
\end{lemma}

\subsection{Separatrices and first integrals}\label{sec:fi}

A \emph{separatrix} of $\Fcal$ is an irreducible germ of analytic curve
$\Gamma\ni0$ which is invariant by $\Fcal$.

\begin{lemma}\label{lem:sep}
Let $\Gamma$ be a separatrix of $\Fcal$.
\begin{enumerate}[label=\textup{(\alph*)}]
\item Near a point $x\in\Gamma\setminus\{0\}$, $\Gamma$ is the plaque
through $x$ of a flow box.
\item Every neighbourhood of $0$ contains a representative $\Fcal_U$ such
that $\Gamma^\circ_U=(\Gamma\cap U)\setminus\{0\}$ is a leaf of $\Fcal_U$.
\end{enumerate}
\end{lemma}

\begin{proof}
(a) A flow box $f$ at $x$ is constant on the regular part of each local
branch of the germ of $\Gamma$ at $x$, since that part is connected and
invariant. Hence each local branch is contained in, and therefore equal
to, the smooth curve $\{f=f(x)\}$, so that the germ of $\Gamma$ at $x$ is
this curve. (b) For $U$ small, $\Gamma\cap U$ is
irreducible, so that $\Gamma^\circ_U$ is connected \cite{Chi89}; by (a) it
is open in the leaf containing it, and it is closed in that leaf because
$\Gamma\cap U$ is closed in $U$.
\end{proof}

Such representatives are \emph{adapted to $\Gamma$}.

\begin{lemma}\label{lem:fi}
Let $\Fcal_U$ be a representative and $H$ a non-constant holomorphic
first integral on $U$.
\begin{enumerate}[label=\textup{(\alph*)}]
\item $H$ is constant on the leaves of $\Fcal_U$, and the separatrices of
$\Fcal$ are the irreducible components of the germ of $\{H=H(0)\}$; there
are finitely many, and at least one.
\item Let $\Gamma$ be a separatrix and $\Sigma$ a transverse section, with
coordinate $z$, at a point of $\Gamma^\circ_{U'}$, where $U'\subset U$ is
adapted to $\Gamma$. Then $H|_\Sigma=H(0)+c\,z^m(1+O(z))$ with $c\neq0$,
and $m=:\ord_\Gamma H\geq1$ does not depend on the choices, including
that of $U'$: two representatives adapted to $\Gamma$ share regular points
of $\Gamma$ arbitrarily close to $0$.
\item If $f_i$ is a reduced equation of the separatrix $\Gamma_i$, then
$H-H(0)=v\prod_if_i^{m_i}$ with $v$ a unit and $m_i=\ord_{\Gamma_i}H$.
\item If $d\geq1$ divides every $\ord_\Gamma H$, there is a holomorphic
first integral $K$ with $K^d=H-H(0)$.
\item Near $0$, $\sing(\Fcal)\subset\{dH=0\}\subset\{H=H(0)\}$.
\end{enumerate}
\end{lemma}

\begin{proof}
(a) In a flow box, $H=h\circ f$ with $h$ holomorphic, so $H$ is constant
on plaques, hence on leaves. A component $Y$ of $\{H=H(0)\}$ is, near a
generic point, a union of plaques, since the zeros of $h-H(0)$ are
isolated; so $Y$ is invariant. Conversely $H$ is constant on the leaf
$\Gamma^\circ_{U'}$, which accumulates at $0$, so $\Gamma\subset\{H=H(0)\}$.
(b) Near a point of $\Gamma^\circ_{U'}$, $H=h\circ f$ with
$h(t)-H(0)=t^m\cdot(\text{unit})$, and $f|_\Sigma$ is a coordinate; the
integer $m$ is locally constant along the connected leaf
$\Gamma^\circ_{U'}$. (c) By (a) the irreducible factors of $H-H(0)$ are the
$f_i$; at a generic point of $\Gamma_i$, where $f_i|_\Sigma$ has a simple
zero, (b) identifies the exponent. (d) Take $K=v^{1/d}\prod f_i^{m_i/d}$;
then $dK\wedge\omega=0$ off $\{K=0\}$, hence everywhere. (e) Where
$dH\neq0$ we have $\omega=a\,dH$ with $a$ holomorphic, and $a$ does not
vanish since $\omega$ has an isolated zero; and $H$ is constant on each
irreducible component of $\{dH=0\}$ through $0$.
\end{proof}

A first integral $H$ is \emph{normalized} if $H(0)=0$; replacing $H$ by
$H-H(0)$ we may always assume this, and then (d) reads $H=K^d$. A
non-constant first integral $H$ is \emph{not a power}, in the
terminology of Mattei and Moussu \cite[Introduction, \S2]{MM80}, if the
exponents $m_i$ in (c) are coprime; by (d), every non-constant
normalized first integral is a power of one which is not a power. A first integral $F_0$ is \emph{primitive} if every
holomorphic first integral is of the form $h\circ F_0$.

\begin{lemma}[{\cite[Introduction, \S2, and Chapitre~V, Th\'eor\`eme~3]{MM80}}]
\label{lem:primitive}
If $\Fcal$ has a non-constant holomorphic first integral, a first integral
$F_0$ is primitive if and only if $F_0-F_0(0)$ is not a power, and such
first integrals exist. If $F_0$ is primitive and $H$ is a non-constant
first integral, then $\ord_\Gamma H\geq\ord_\Gamma F_0$ for every
separatrix $\Gamma$.
\end{lemma}

The first assertion is the factorization theorem of Mattei and Moussu,
together with Lemma~\ref{lem:fi}(d); the last follows from
$\ord_\Gamma(h\circ F_0)=\ord_{F_0(0)}(h-h(F_0(0)))\cdot\ord_\Gamma F_0$.

\subsection{Sectors, leaf counts and leading terms}\label{sec:sectors}

Sectors are defined in Definition~\ref{def:sector}; the same definition
applies to a coordinate $z$ centred at $0$ in $\C$, the sector then having
its vertex at $0$. As in \cite{Was65,RS21}, a \emph{proper subsector} of
$S(r;\theta_1,\theta_2)$ is a set
$\{0<|z|\leq r',\ \arg z\in[\psi_1,\psi_2]+2\pi\Z\}$ with $0<r'<r$ and
$\theta_1<\psi_1\leq\psi_2<\theta_2$.

Let $\sim$ be an equivalence relation on a neighbourhood of $0$ in $\C$ and
$S=S(r;\theta_1,\theta_2)$ a sector contained in it. We put
$\nu(S,\rho)=\sup_C\#(C\cap S\cap\{|z|>\rho\})$ and
$\nu(S)=\sup_C\#(C\cap S)$, where $C$ ranges over the equivalence classes.
For a representative $\Fcal_U$ and a sector in a transverse section, these
are the numbers $\nu_U(S,\rho)$ and $\nu_U(S)$ of
Definition~\ref{def:counts}, relative to the relation $\sim_U$; for a
foliation $\Gcal$ on a complex surface $X$ we write $\nu_X(S,\rho)$ and
$\nu_X(S)$. As in Definition~\ref{def:counts}, $S$ has \emph{sub-logarithmic
counts} if $\liminf_{\rho\to0}\nu(S,\rho)/\log(1/\rho)=0$. Bounded counts
are sub-logarithmic.

\begin{lemma}\label{lem:counts}
\begin{enumerate}[label=\textup{(\alph*)}]
\item If $S'\subset S$, then $\nu(S',\rho)\leq\nu(S,\rho)$ for every
$\rho$.
\item If $U'\subset U$ and $S$ is a sector in a transverse section
contained in $U'$, then $\nu_{U'}(S,\rho)\leq\nu_U(S,\rho)$ for every
$\rho$.
\item Let $S$ be a sector in a transverse section $\Sigma$, with coordinate
$z$. Let $S'$ be a subset of a transverse section $\Sigma'$, with
coordinate $w$, and $\beta\colon S'\to S$ an injective map such that
$\beta(w)\sim_Uw$ and $|z(\beta(w))|\geq c|w|^\kappa$ for all $w\in S'$,
where $c,\kappa>0$. Then, for every $\rho>0$, every leaf of $\Fcal_U$
meets $S'\cap\{|w|>\rho\}$ in at most $\nu_U(S,c\rho^\kappa)$ points. In
particular, every leaf meets $S'$ in at most $\nu_U(S)$ points, and if
$S'$ is a sector and $S$ has sub-logarithmic counts, then so has $S'$.
\item In the situation of \textup{(c)}, let $S$ have opening $\vartheta$
and let $\sigma\colon(\Sigma',a')\to(\Sigma,a)$ be a germ of biholomorphism
with $\sigma(w)\sim_U w$ for $w$ close to $a'$, for instance a holonomy
map. Then, for any coordinate on $\Sigma'$, the section $\Sigma'$ contains
sectors $S'$ of every opening $\vartheta'<\vartheta$ such that
$\nu_U(S',\rho)\leq\nu_U(S,c\rho)$ for some $c>0$ and every $\rho$; in the
coordinate $z\circ\sigma$ it contains such a sector of opening
$\vartheta$. In particular, if $S$ has bounded, or sub-logarithmic,
counts, so has $S'$.
\end{enumerate}
\end{lemma}

\begin{proof}
(a) and (b) are clear, since the leaves of $\Fcal_{U'}$ lie in leaves of
$\Fcal_U$. (c) The map $\beta$ sends $L\cap S'\cap\{|w|>\rho\}$
injectively into $L\cap S\cap\{|z|>c\rho^\kappa\}$ for every leaf $L$.
For the last assertion, let $\rho_k\to0$ with
$\nu_U(S,\rho_k)/\log(1/\rho_k)\to0$ and put
$\rho'_k=(\rho_k/c)^{1/\kappa}$. Then
$\nu_U(S',\rho'_k)\leq\nu_U(S,\rho_k)$, while
$\log(1/\rho'_k)=\kappa^{-1}\bigl(\log(1/\rho_k)+\log c\bigr)$, so that
$\nu_U(S',\rho'_k)/\log(1/\rho'_k)\to0$. (d) Let
$S=S(r;\theta_1,\theta_2)$ in the coordinate $z$, let $w$ be a coordinate
on $\Sigma'$, and write $z(\sigma(w))=\lambda w(1+O(w))$ with
$\lambda\neq0$. Given $\vartheta'<\vartheta$, put
$\delta=(\vartheta-\vartheta')/3$ and
$S'=S(r';\theta_1-\arg\lambda+\delta,\theta_1-\arg\lambda+\delta+\vartheta')$.
Since $\arg z(\sigma(w))=\arg\lambda+\arg w+O(|w|)$ and
$|z(\sigma(w))|=|\lambda||w|(1+O(|w|))$, we have $\sigma(S')\subset S$
and $|z(\sigma(w))|\geq|\lambda||w|/2$ on $S'$ for $r'$ small, and
$\sigma$ is injective on $S'$. Apply (c) with $\beta=\sigma|_{S'}$,
$c=|\lambda|/2$ and $\kappa=1$. In the coordinate $w=z\circ\sigma$ the map
$\sigma$ is the identity, and $S'=S(r';\theta_1,\theta_2)$ will do.
\end{proof}

We say that a holomorphic function $u$ on a sector $S$ has a \emph{leading
term} if $u(z)=b+c\,z^\alpha(1+\varepsilon(z))$, where $c\neq0$,
$\Re\alpha>0$, $z^\alpha$ is computed with a branch of $\log z$ on $S$, and
$\varepsilon(z)\to0$ as $z\to0$, uniformly on every proper subsector. A
moderate sectorial first integral (Corollary~\ref{cor:modA}) has a
leading term $a_dz^d$, where $d\geq1$ is the smallest positive index with
$a_d\neq0$: indeed
$|\varphi(z)-a_0-a_dz^d|\leq A_{d+1}|z|^{d+1}$ on each proper subsector.

\begin{lemma}\label{lem:leading}
Let $u$ be holomorphic on a sector $S=S(r;\theta_1,\theta_2)$ and have a
leading term $c\,z^\alpha$. Let $[\psi_1,\psi_2]\subset(\theta_1,\theta_2)$
with $\psi_2-\psi_1<2\pi\Re\alpha/|\alpha|^2$. Then $u$ is injective on
$S(r';\psi_1,\psi_2)$ for some $r'>0$. Consequently, if $u$ is a sectorial
first integral of $\Fcal_U$, every leaf of $\Fcal_U$ meets
$S(r';\psi_1,\psi_2)$ in at most one point.
\end{lemma}

\begin{proof}
Choose $[\psi'_1,\psi'_2]\subset(\theta_1,\theta_2)$ containing
$[\psi_1,\psi_2]$ in its interior, and $\kappa\in(0,1/2)$ such that, for
$z$ in $S(r/4;\psi_1,\psi_2)$, the closed disc of radius $\kappa|z|$ about
$z$ lies in the proper subsector
$S_1=\{0<|z|\leq r/2,\ \arg z\in[\psi'_1,\psi'_2]+2\pi\Z\}$ of $S$. By the
Cauchy inequality, $|z\,\varepsilon'(z)|\leq\kappa^{-1}\max\{|\varepsilon(w)|:
|w-z|=\kappa|z|\}$, which tends to $0$ as $z\to0$ in
$S(r/4;\psi_1,\psi_2)$, because $\varepsilon\to0$ uniformly on $S_1$.

Put $z=e^\zeta$, with $\zeta$ in the convex half-strip
$P=\{\Re\zeta<\log r',\ \psi_1<\Im\zeta<\psi_2\}$. For $r'$ small,
$u(e^\zeta)-b=\exp\Lambda(\zeta)$ with
$\Lambda(\zeta)=\log c+\alpha\zeta+\log(1+\varepsilon(e^\zeta))$, and
$\Lambda'=\alpha(1+\eta)$ with $\eta_0(r'):=\sup_P|\eta|\to0$ as
$r'\to0$. Integrating along segments,
$\Lambda(\zeta_1)-\Lambda(\zeta_2)=\alpha(1+\bar\eta)(\zeta_1-\zeta_2)$ with
$|\bar\eta|\leq\eta_0(r')$. If $u(e^{\zeta_1})=u(e^{\zeta_2})$, then
$\Lambda(\zeta_1)-\Lambda(\zeta_2)=2\pi ik$, $k\in\Z$. If $k\neq0$, then
$|\Im(\zeta_1-\zeta_2)|\geq2\pi\Re\bigl(1/(\alpha(1+\bar\eta))\bigr)$,
which tends to $2\pi\Re\alpha/|\alpha|^2>\psi_2-\psi_1$ as $r'\to0$,
uniformly; this contradicts $|\Im(\zeta_1-\zeta_2)|<\psi_2-\psi_1$ for $r'$
small. Hence $k=0$ and $\zeta_1=\zeta_2$, and $u$ is injective on
$S(r';\psi_1,\psi_2)$. The last assertion follows, as a sectorial first
integral is constant on the trace of each leaf.
\end{proof}

\subsection{Finite groups of germs of diffeomorphisms}
\label{sec:onevar-lemmas}

The next lemma consists of Propositions~1.1 and~1.2 of Mattei and Moussu
\cite[Chapitre~I, \S1, pp.~475--476]{MM80}.

\begin{lemma}\label{lem:finite}
\begin{enumerate}[label=\textup{(\alph*)}]
\item Let $G$ be a finite subgroup of $\Diff(\C,0)$, of order $n$. Then
there is a coordinate $\xi$ in which $G$ is the group of the rotations
$\xi\mapsto\zeta\xi$, with $\zeta^n=1$. Hence the map $g\mapsto g'(0)$ is
an isomorphism of $G$ onto the group $\mu_n$ of the $n$-th roots of unity.
\item Let $Q\in\Ocal_1$ be non-constant, of order
$d=\ord_0(Q-Q(0))$, and let $w$ be a coordinate with $Q=Q(0)+w^d$. Then
the invariance group
$\Inv(Q)=\{\gamma\in\Diff(\C,0):Q\circ\gamma=Q\}$ is the cyclic group of
order $d$ formed by the germs $w^{-1}(\zeta w)$, $\zeta^d=1$.
\end{enumerate}
In particular, a germ of finite order is analytically linearizable, its
order equals the order of its multiplier in $\C^*$, and a germ of finite
order tangent to the identity is the identity.
\end{lemma}

\begin{lemma}\label{lem:arc}
Let $m\geq1$, $\vartheta\in(0,2\pi]$ and
$I=(\theta_1,\theta_1+\vartheta)$.
\begin{enumerate}[label=\textup{(\alph*)}]
\item For every $\psi\in\R$, the number of indices $j\in\{0,\dots,m-1\}$
with $\psi+2\pi j/m\in I+2\pi\Z$ is at most
$\lceil m\vartheta/2\pi\rceil$.
\item There are $\psi\in\R$ and $\delta>0$ such that at least
$\lceil m\vartheta/2\pi\rceil$ indices $j\in\{0,\dots,m-1\}$ satisfy
$\psi+2\pi j/m\in(\theta_1+\delta,\theta_1+\vartheta-\delta)+2\pi\Z$.
\end{enumerate}
\end{lemma}

\begin{proof}
(a) Let $A=\psi+(2\pi/m)\Z$. Since $I$ has length at most $2\pi$, each
class modulo $2\pi$ has at most one representative in $I$, and the classes
of $\psi+2\pi j/m$, $0\leq j<m$, are distinct. Hence the number in question
is $\#(A\cap I)$. If $A\cap I$ has $k\geq1$ elements, they are consecutive
terms of an arithmetic progression of step $2\pi/m$ contained in an open
interval of length $\vartheta$, so $2\pi(k-1)/m<\vartheta$, that is,
$k-1<m\vartheta/2\pi$. Since $k-1$ is an integer,
$k\leq\lceil m\vartheta/2\pi\rceil$. (b) Let
$k=\lceil m\vartheta/2\pi\rceil$. Then $k\leq m$, because
$\vartheta\leq2\pi$, and $2\pi(k-1)/m<\vartheta$. Put
$\delta=\bigl(\vartheta-2\pi(k-1)/m\bigr)/4>0$ and $\psi=\theta_1+2\delta$.
For $0\leq j\leq k-1$ the numbers $\psi+2\pi j/m$ lie in
$[\theta_1+2\delta,\theta_1+\vartheta-2\delta]$.
\end{proof}

\begin{lemma}\label{lem:periodic}
Let $g\in\Diff(\C,0)$ have finite order $m$, and let
$S=S(r;\theta_1,\theta_2)$ be a sector of opening $\vartheta$. For every
$\epsilon>0$ there is $z_0\in S$ such that the orbit
$\{g^j(z_0)\}_{0\leq j<m}$ consists of $m$ distinct points of modulus less
than $\epsilon$, and at least $\lceil m\vartheta/2\pi\rceil$ of them lie in
$S$.
\end{lemma}

\begin{proof}
Let $\xi$ be the coordinate of Lemma~\ref{lem:finite}(a) for the group
generated by $g$, normalized by $\xi'(0)=1$; in it $g$ is the rotation
by the primitive $m$-th root of unity $\lambda=g'(0)$. The orbit of
$z_0=\xi^{-1}(te^{i\psi})$, $t>0$ small, is therefore the set of the $m$
distinct points $\xi^{-1}(te^{i(\psi+2\pi j/m)})$, $0\leq j<m$. Choose
$\psi$ and $\delta$ as in Lemma~\ref{lem:arc}(b) for the interval
$(\theta_1,\theta_2)$. Since $\xi^{-1}(w)=w(1+O(w))$, there is
$t_0>0$ such that $|\xi^{-1}(w)|<2|w|$ and $|\arg(\xi^{-1}(w)/w)|<\delta$
for $0<|w|<t_0$. If $t<\min\{t_0,\epsilon/2,r/2\}$, all points of the orbit
have modulus less than $\min\{\epsilon,r\}$, and those with
$\psi+2\pi j/m\in(\theta_1+\delta,\theta_2-\delta)+2\pi\Z$ lie in $S$.
\end{proof}
\subsection{Reduction of singularities in dimension two}\label{sec:reduction}

Let $n=2$ and $\omega=a\,dx+b\,dy$. The foliation is also generated by the
vector field $X=b\,\partial_x-a\,\partial_y$, and
$dH\wedge\omega=X(H)\,dx\wedge dy$. The singular point $0$ is \emph{reduced}
(\emph{simple} in \cite{RS21}, \emph{irreducible} in \cite{Sca16}) if
the eigenvalues $\lambda_1,\lambda_2$ of the linear part of $X$ satisfy one
of the following conditions: $\lambda_1\lambda_2\neq0$ and
$\lambda_1/\lambda_2\notin\Q_{>0}$ (a \emph{non-degenerate} singularity),
or exactly one of them vanishes (a \emph{saddle-node}). The
\emph{irreducible} forms of Mattei and Moussu
\cite[Chapitre~II, p.~480]{MM80} are defined by the weaker condition
$\lambda_1/\lambda_2,\lambda_2/\lambda_1\notin\N$ in the first case, so
that their results on irreducible forms apply to reduced singularities.

By Seidenberg's theorem \cite{Sei68} there is a finite composition of
blow-ups of points $\pi\colon(M,D)\to(\C^2,0)$, $D=\pi^{-1}(0)$, such that
the strict transform $\tF$, defined locally by the saturation of
$\pi^*\omega$, has only reduced singularities; we call $\pi$ a
\emph{reduction of singularities} of $\Fcal$. The \emph{minimal} reduction
is obtained by blowing up, at each step, the singular points which are not
reduced. If $0$ is not reduced, the exceptional divisor $D$ is a connected
curve with normal crossings whose irreducible components are smooth
rational curves, any two of them meeting in at most one point, called a
\emph{corner}; its dual graph is a tree. Each component of $D$ is either
$\tF$-invariant or \emph{dicritical}, that is, generically transverse to
$\tF$. The germ $\Fcal$ is \emph{non-dicritical} if all the components of
$D$ are invariant, and a \emph{generalized curve} if no singular point of
$\tF$ is a saddle-node \cite{CLS84}; when $0$ is reduced, $\Fcal$ is
non-dicritical, and it is a generalized curve if and only if $0$ is not a
saddle-node. These properties do not depend on the reduction \cite{CLS84}.

We use the following classical facts.
\begin{enumerate}[label=\textup{(R\arabic*)}]
\item\label{R1} \emph{Non-degenerate reduced singularities}
\cite[Chapitre~II, \S1, and Appendice~II, Th\'eor\`eme~1]{MM80}. There are
holomorphic coordinates in which $X=xA\,\partial_x+yB\,\partial_y$ with
$A(0)=\lambda_1$ and $B(0)=\lambda_2$. Thus the axes are separatrices; they
are smooth, transverse, and tangent to the eigendirections. There are
exactly two separatrices: an invariant curve parametrized by
$t\mapsto(c\,t^a+\cdots,d\,t^b+\cdots)$, with $cd\neq0$ and $a,b\geq1$,
would give $\lambda_2/\lambda_1=b/a\in\Q_{>0}$. The holonomy map of
$\{y=0\}$ along a positive loop around $0$, computed on $\{x=x_0\}$ with
the coordinate $y$, has multiplier $\exp(2\pi i\lambda_2/\lambda_1)$
\cite[p.~482]{MM80}.
\item\label{R2} \emph{Saddle-nodes}
\cite[Chapitre~II, \S1 and \S3, and Appendice~II, Th\'eor\`eme~2]{MM80}. If
$0$ is a saddle-node, there is a unique separatrix tangent to the
eigendirection of the nonzero eigenvalue, the \emph{strong separatrix}.
The other formal invariant curve, tangent to the kernel of the linear
part, is the \emph{weak separatrix}; it need not converge \cite{MR82}.
There are coordinates in which $\Fcal$ is defined by
$y^m\,dx-(x+yA(x,y))\,dy$, with $m\geq2$
\cite[formula~(1.2), p.~480]{MM80}. The strong separatrix is then
$\{y=0\}$, and its holonomy map along the loop
$\theta\mapsto(e^{2\pi i\theta},0)$, computed on $\{x=1\}$ with the
coordinate $y$, is
\[
h(y)=y+2\pi i\,y^m+O(y^{m+1})
\]
\cite[p.~482]{MM80}; see also \cite[\S4.3]{FMR24}. Thus $h$ is tangent to
the identity and different from the identity, and it has infinite order
(Lemma~\ref{lem:finite}). Consequently, a saddle-node has no non-constant
holomorphic first integral
\cite[Chapitre~II, Th\'eor\`eme~$B_0$, and \S3, p.~485]{MM80}: the
restriction of such a first integral to a transverse section to the strong
separatrix would be a non-constant germ invariant under $h$, contradicting
Lemma~\ref{lem:finite}(b).
\item\label{R3} \emph{Linearization} \cite[Chapitre~II, \S2,
Th\'eor\`eme~2, p.~482]{MM80}; see also \cite[p.~595]{MR83} and
\cite[Theorem~3.1]{RS21}. Let $0$ be a non-degenerate reduced singularity
such that the holonomy of one of its separatrices is periodic. By
\ref{R1}, the corresponding ratio of eigenvalues is rational, since the
multiplier of the holonomy is a root of unity, and it is negative because
the singularity is reduced. A periodic germ is analytically linearizable,
so the cited theorem gives holomorphic coordinates in which $\Fcal$ is
defined by $p\,y\,dx+q\,x\,dy$, with $p,q\geq1$ coprime.
\end{enumerate}

These notions and facts apply to the germ, at a singular point $s$, of a
holomorphic foliation $\Gcal$ with isolated singularities on a complex
surface $X$. The next lemma rests on the theorem of Mattei and Moussu on
the connectedness of the fibres of a function which is not a power
\cite[Chapitre~V, Lemme~1, p.~500]{MM80}. Its part~(b) is the form in which
this theorem is used in the proof of Th\'eor\`eme~1 of
\cite[Chapitre~V, \S3, p.~504]{MM80}.

\begin{lemma}[Connected fibres]\label{lem:connected}
Let $\Gcal$ be a holomorphic foliation with isolated singularities on a
complex surface $X$, let $s\in X$, and let $g$ be a holomorphic first
integral of\/ $\Gcal$ on a neighbourhood of $s$, with $g(s)=0$, whose germ
at $s$ is not a power. Then every neighbourhood of $s$ contains an open
neighbourhood $V$ of $s$ with the following properties.
\begin{enumerate}[label=\textup{(\alph*)}]
\item For every $c\in\C^*$, the set $\{g=c\}\cap V$ is connected
(possibly empty) and contained in a leaf of\/ $\Gcal$.
\item If $x\in V\setminus\sing(\Gcal)$, $g(x)=0$ and $\tau$ is a
transverse section at $x$, then
$\Inv(g|_\tau)\subset\Holv(\Gcal,\tau,x)$.
\end{enumerate}
\end{lemma}

\begin{proof}
Let $B$ be a connected open neighbourhood of $s$, contained in the given
neighbourhood, on which $g$ is a first integral, such that $\Gcal$ has no
singular point in $B\setminus\{s\}$ and the critical points of $g$ in $B$
lie in $\{g=0\}$ (Lemma~\ref{lem:fi}(e)). By the lemma of Mattei and
Moussu quoted above, $B$ contains an open neighbourhood $V$ of $s$ such
that the level sets of $g|_V$ are connected.

(a) Let $c\neq0$. The set $Y=\{g=c\}\cap V$ is a connected smooth curve
contained in $V\setminus\{s\}$, and it is tangent to $\Gcal$, because
$dg\wedge\omega=0$ for a local generator $\omega$ of\/ $\Gcal$. Hence every
flow box submersion of\/ $\Gcal$ is locally constant on $Y$, so that the
set of the points of $Y$ which lie on a given leaf is open in $Y$. Since
$Y$ is connected, it is contained in one leaf.

(b) In a flow box $f$ at $x$ we have $g=h\circ f$, where $h$ is not
constant because $B$ is connected; as $f|_\tau$ is a coordinate, $x$ is an
isolated zero of $g|_\tau$. Let $\gamma\in\Inv(g|_\tau)$. For $w\in\tau$
close to $x$, $w\neq x$, the points $w$ and $\gamma(w)$ belong to
$\{g=c\}\cap V$ with $c=g(w)\neq0$, and by (a) they lie on the same leaf.
Hence $\gamma\in\Holv(\Gcal,\tau,x)$.
\end{proof}

\subsection{Germs with a first integral}\label{sec:sections}

\begin{lemma}\label{lem:integrable-gc}
A germ of singular holomorphic foliation at $0\in\C^2$ with a non-constant
holomorphic first integral $H$ is a non-dicritical generalized curve.
\end{lemma}

\begin{proof}
Let $\pi$ be a reduction of singularities, with strict transform $\tF$.
Since $\pi$ is a biholomorphism off the exceptional divisor $D$, the
function $H\circ\pi$ is a non-constant first integral of $\tF$ near $D$,
and it equals $H(0)$ on $D$. If a component of $D$ were dicritical, at a
generic point of it a flow box $f$ would restrict to a coordinate on the
component, and $H\circ\pi=h\circ f$ would force $h$, hence $H\circ\pi$, to
be constant near that point, and then everywhere. A saddle-node of $\tF$
would carry the non-constant first integral $H\circ\pi$, contradicting
\ref{R2}.
\end{proof}

\section{Invariance groups in one variable: proof of Theorem~\ref{thm:onevar}}\label{sec:onevar}

Transverse sections to a foliation of codimension one are one-dimensional,
and the leaves of $\Fcal_U$ induce on each of them an equivalence relation,
whose invariance group is the virtual holonomy group
(Definition~\ref{def:G}). Throughout this section, $\sim$ is an
equivalence relation on an open neighbourhood $\Delta$ of $0$ in $\C$,
$G=\Inv(\sim)$, all sectors have vertex $0$ and are contained in $\Delta$,
and $\nu(S,\rho)$, $\nu(S)$ are as in Section~\ref{sec:sectors}. If
$g\in G$, there is a neighbourhood $V$ of $0$ with $g(z)\sim z$ for
$z\in V$; consequently, for every finite set of elements of $G$ there is a
neighbourhood of $0$ on which all of them move the points within their
classes.

\subsection{Bounded counts}

\begin{proof}[Proof of Theorem~\textup{\ref{thm:onevar}(a)}]
Let $S=S(r;\theta_1,\theta_2)$ be a sector of opening $\vartheta$ with
$\nu(S)\leq N$.

\emph{Step 1: every element of $G$ has finite order.} Let $g\in G$ and
$\lambda=g'(0)$. By Dirichlet's approximation theorem there is an integer
$m\geq1$ such that the argument of $\lambda^m$ in $(-\pi,\pi]$ satisfies
$|\arg\lambda^m|<\vartheta/(4N)$. Put $f=g^m\in G$ and
$\mu=f'(0)=\lambda^m$. Let $\theta_c=(\theta_1+\theta_2)/2$ and
$S_0=S(r_0;\theta_c-\vartheta/8,\theta_c+\vartheta/8)$. For
$0\leq j\leq N$ we have $f^j(z)=\mu^jz(1+O(z))$, hence
\[
\arg f^j(z)\in\arg z+j\arg\mu+O(|z|)+2\pi\Z,\qquad
|f^j(z)|=|\mu|^j|z|\,(1+O(|z|)),
\]
with $|j\arg\mu|<\vartheta/4$. If $r_0$ is small enough, then for every
$z\in S_0$ and $0\leq j\leq N$: the point $f^j(z)$ belongs to $S$, because
its argument differs from $\theta_c$ by less than
$\vartheta/8+\vartheta/4+\vartheta/8=\vartheta/2$ and its modulus is less
than $r$; $f^j(z)\sim z$, because $f^j\in G$; and $f^i$ is injective on a
disc containing all the points $f^l(z)$, $0\leq l\leq N$. Since the class
of $z$ meets $S$ in at most $N$ points, there are $0\leq i<j\leq N$ with
$f^i(z)=f^j(z)=f^i(f^{j-i}(z))$, so $f^{j-i}(z)=z$. Hence
\[
S_0\subset\bigcup_{d=1}^{N}\Fix(f^d),\qquad\Fix(f^d)=\{z:\ f^d(z)=z\}.
\]
If no $f^d$ with $1\leq d\leq N$ were the identity, each $\Fix(f^d)$ would
be a discrete subset of a disc centred at $0$, being the zero set of the
nonzero holomorphic function $f^d(z)-z$; a finite union of discrete sets is
countable and cannot contain the open set $S_0$. Hence $f^d=\id$ for some
$d\leq N$, and $g^{md}=\id$.

\emph{Step 2: if $g\in G$ has finite order $m$, then
$\lceil m\vartheta/2\pi\rceil\leq N$.} Choose $\epsilon>0$ such that
$g^j(z)\sim z$ for $|z|<\epsilon$ and $0\leq j<m$. By
Lemma~\ref{lem:periodic} there is $z_0\in S$ whose orbit consists of points
of modulus less than $\epsilon$, at least $\lceil m\vartheta/2\pi\rceil$ of
which lie in $S$. They all belong to the class of $z_0$. Hence
$\lceil m\vartheta/2\pi\rceil\leq N$, that is, $m\vartheta/2\pi\leq N$.

\emph{Step 3: $G$ is finite cyclic of order at most $2\pi N/\vartheta$.}
By Steps~1 and~2, every $g\in G$ has finite order at most
$M=\lfloor2\pi N/\vartheta\rfloor$. By Lemma~\ref{lem:finite}, $g'(0)$ is a
root of unity of the same order as $g$, and $g'(0)=1$ only if $g=\id$.
Hence $g\mapsto g'(0)$ is an injective homomorphism from $G$ into the group
$\mu_{M!}$ of $M!$-th roots of unity. Thus $G$ is finite and cyclic, and
Step~2 applied to a generator gives $|G|\leq2\pi N/\vartheta$.

\emph{Sharpness.} Let $N\geq1$, $\vartheta\in(0,2\pi]$, and
$m=\lfloor2\pi N/\vartheta\rfloor\geq1$. Define $z\sim_mw$ if $z^m=w^m$. If
$\gamma\in\Inv(\sim_m)$, then $(\gamma(z)/z)^m=1$ near $0$; the function
$\gamma(z)/z$ extends continuously to $0$, so it is a constant $m$-th root
of unity. Hence $\Inv(\sim_m)$ is the group of rotations
$z\mapsto\zeta z$, $\zeta^m=1$, of order $m$. Every class other than
$\{0\}$ consists of $m$ points equally spaced on a circle, so by
Lemma~\ref{lem:arc}(a) every sector of opening $\vartheta$ contains at most
$\lceil m\vartheta/2\pi\rceil\leq N$ of them.
\end{proof}

Step~2 of the proof does not use the hypothesis of (a), and we record it
for later use: \emph{if $\Inv(\sim)$ contains an element of order $m$,
then every sector $S$ of opening $\vartheta$ satisfies
$\nu(S)\geq\lceil m\vartheta/2\pi\rceil$.}

\subsection{Sub-logarithmic counts}

The proof of Theorem~\ref{thm:onevar}(b) treats separately the elements
whose multiplier has modulus different from one, which produce classes
with a number of points of order $\log(1/\rho)$, and those whose
multiplier has modulus one, for which Step~1 above is made quantitative.

\begin{lemma}\label{lem:hyperbolic}
Let $S$ be a sector of opening $\vartheta$ and $g\in G$ with
$|g'(0)|\neq1$. Then there are $c>0$ and $\rho_0>0$ such that
$\nu(S,\rho)\geq c\log(1/\rho)$ for $0<\rho<\rho_0$.
\end{lemma}

\begin{proof}
Replacing $g$ by $g^{-1}$, we may assume that $\lambda=g'(0)$ satisfies
$0<|\lambda|<1$. By the theorem of Koenigs there is a germ of coordinate
$\phi(z)=z+O(z^2)$ with $\phi\circ g=\lambda\phi$. Let
$S=S(r;\theta_1,\theta_2)$ and $\theta_c=(\theta_1+\theta_2)/2$. Choose
$r_1\in(0,r)$ such that $g(z)\sim z$ for $|z|<r_1$, $\phi$ is injective on
$\{|z|<r_1\}$, its image contains $\{|w|<r_1/2\}$, and
$|z|/2\leq|\phi(z)|\leq2|z|$ and $|\arg(\phi(z)/z)|<\vartheta/8$ for
$0<|z|<r_1$. Fix $w_0$ with $|w_0|=r_1/8$ and $\arg w_0=\theta_c$, and put
$z_j=\phi^{-1}(\lambda^jw_0)$ for $j\geq0$. Then $z_{j+1}=g(z_j)$ and
$|z_j|<r_1$, so that all the $z_j$ belong to the class of $z_0$. The point
$z_j$ lies in $S$ as soon as
$\theta_c+j\arg\lambda\in(\theta_c-\vartheta/4,\theta_c+\vartheta/4)+2\pi\Z$,
and $|z_j|\geq|\lambda|^j|w_0|/2>\rho$ as soon as $j\leq J(\rho)$, where
$J(\rho)$ is the integer part of $\log(r_1/(16\rho))/\log(1/|\lambda|)$.
If $\arg\lambda/2\pi$ is a rational number with denominator $q$, the
first condition holds for every multiple $j$ of $q$; otherwise, by the
equidistribution theorem, it holds for at least
$(\vartheta/8\pi)J$ of the indices $0\leq j\leq J$, when $J$ is large.
In both cases, the class of $z_0$ has at least $c'J(\rho)$ points in
$S\cap\{|z|>\rho\}$ for $\rho$ small, with $c'>0$ independent of
$\rho$, and $J(\rho)\geq\log(1/\rho)/(2\log(1/|\lambda|))$ for $\rho$
small.
\end{proof}

\begin{lemma}\label{lem:elliptic}
Let $S$ be a sector of opening $\vartheta$ and $g\in G$ with
$|g'(0)|=1$. There is $K>0$, depending on $g$ and $S$, with the following
property: if $\rho>0$ is so small that $M=\nu(S,\rho)$ satisfies
$1\leq M<\infty$ and $M^2\rho<K$, then $g$ has finite order.
\end{lemma}

\begin{proof}
Let $S=S(r;\theta_1,\theta_2)$, $\theta_c=(\theta_1+\theta_2)/2$, and
$\lambda=g'(0)$. Choose $r_1\in(0,r/2)$ and $C\geq1/r_1$ such that
$g(z)\sim z$ for $|z|<r_1$, $g$ is injective on $\{|z|<r_1\}$, and
$g(z)=\lambda z\bigl(1+h(z)\bigr)$ with $|h(z)|\leq C|z|$ for $|z|<r_1$.

\emph{Iteration estimate.} Let $k\geq1$ and $|z|\leq1/(8kC)$. Then, for
$0\leq i\leq k$, the points $g^i(z)$ are defined, belong to the class of
$z$, and satisfy
\[
e^{-1}|z|\leq|g^i(z)|\leq e^{1/2}|z|,\qquad
\bigl|\arg g^i(z)-\arg z-i\arg\lambda\bigr|\leq8Ck|z|\pmod{2\pi}.
\]
Indeed, by induction on $i$, if $|g^l(z)|\leq2|z|$ for $l<i$, then
$|h(g^l(z))|\leq2C|z|\leq1/4$, so that
$\bigl|\log|1+h(g^l(z))|\bigr|\leq4C|z|$ and
$|\arg(1+h(g^l(z)))|\leq4C|z|$; summing over $l<i\leq k$ gives
$\bigl|\log(|g^i(z)|/|z|)\bigr|\leq4Ck|z|\leq1/2$ and the estimate on the
argument, and in particular $|g^i(z)|\leq2|z|<r_1$.

\emph{Choice of a power.} Let $M=\nu(S,\rho)$ be as in the statement and
$Q=\lceil16\pi M/\vartheta\rceil\leq18\pi M/\vartheta$. By Dirichlet's
theorem there is $1\leq m\leq Q$ such that $m\arg\lambda\equiv\alpha
\pmod{2\pi}$ with $|\alpha|\leq2\pi/Q\leq\vartheta/(8M)$. Put $k=mM$ and
consider
\[
W=\Bigl\{z:\ e\rho<|z|<\frac{\vartheta}{64\,CQM},\
|\arg z-\theta_c|<\frac\vartheta8\Bigr\}.
\]
Since $Q\leq18\pi M/\vartheta$, the set $W$ is open and non-empty as soon
as $e\rho<\vartheta^2/(1152\pi CM^2)$; we put $K=\vartheta^2/(1152\pi eC)$.
Let $z\in W$ and $z_j=g^{mj}(z)$, $0\leq j\leq M$. Since
$|z|\leq\vartheta/(64Ck)\leq1/(8kC)$, the iteration estimate applies with
$i=mj\leq k$: the points $z_j$ belong to the class of $z$, have modulus
in $(\rho,r)$, and their arguments differ from $\theta_c$ by at most
$\vartheta/8+|j\alpha|+8Ck|z|\leq3\vartheta/8$. Hence the $M+1$ points
$z_0,\dots,z_M$ lie in $S\cap\{|w|>\rho\}$ and in one class, which has at
most $M$ points there. So $z_i=z_j$ for some $0\leq i<j\leq M$, and, $g$
being injective on $\{|w|<r_1\}$ and the iterates staying in this disc,
$g^{m(j-i)}(z)=z$. Thus $W\subset\bigcup_{d=1}^{M}\Fix(g^{md})$, and as in
Step~1 of the proof of Theorem~\ref{thm:onevar}(a) some $g^{md}$ is the
identity.
\end{proof}

\begin{proof}[Proof of Theorem~\textup{\ref{thm:onevar}(b)} and~\textup{(c)}]
(b) Let $\rho_k\to0$ with $\nu(S,\rho_k)/\log(1/\rho_k)\to0$. Since the
classes cover $S$, $\nu(S,\rho_k)\geq1$ for $\rho_k$ smaller than the
radius of $S$. Let $g\in G$. By Lemma~\ref{lem:hyperbolic},
$|g'(0)|=1$. For $k$ large, $M_k=\nu(S,\rho_k)$ is finite and
$M_k^2\rho_k\to0$, so Lemma~\ref{lem:elliptic} shows that $g$ has finite
order. By Lemma~\ref{lem:finite}, an element of finite order with
multiplier $1$ is the identity, so the homomorphism $g\mapsto g'(0)$ is
injective on $G$, and its image consists of roots of unity. A finitely
generated subgroup of $G$ is therefore isomorphic to a group generated by
finitely many roots of unity, which is finite and cyclic.

(c) The class of $z\neq0$ is $\{2^kz:\ k\in\Z,\ |2^kz|<1\}$, and its
points of modulus greater than $\rho$ correspond to the integers $k$ with
$\rho<2^k|z|<1$, which are at most $\log_2(1/\rho)+1$. The germ
$z\mapsto z/2$ maps each class into itself.
\end{proof}

The torsion conclusion of Theorem~\ref{thm:onevar}(b) cannot be improved
to finiteness, even under a much slower growth. Following Rib\'on
\cite[\S7, Remark~25]{Rib26}, let $\varphi(\rho)\to\infty$ as $\rho\to0$,
choose radii $r_1>r_2>\dots\to0$ in $(0,1)$ with $\varphi(\rho)\geq4^k$
for $0<\rho<r_k$, and let $\sim$ be the relation on the unit disc
generated by the rotations $\sigma_k(z)=e^{2\pi i/2^k}z$ restricted to
$\{|z|<r_k\}$. The class of a point with $r_{k+1}\leq|z|<r_k$ consists of
the $2^k$ points $e^{2\pi ij/2^k}z$, so that
$\nu(S,\rho)\leq2^k\leq2^{-k}\varphi(\rho)$ for $r_{k+1}\leq\rho<r_k$,
whereas $\Inv(\sim)$ contains every $\sigma_k$ and is infinite; with
$r_k=e^{-4^k}$, $\nu(S,\rho)\leq\max\{1,(\log(1/\rho))^{1/2}\}$.

\subsection{Orbit growth of germs and finitely generated groups}
\label{sec:growth}

The proof of Theorem~\ref{thm:onevar}(b) gives lower bounds for the counts
produced by a single element of infinite order. They depend on the
multiplier, and when the multiplier is a root of unity no approximation
argument is needed, which improves the exponent.

\begin{lemma}\label{lem:parabolic}
Let $S$ be a sector of opening $\vartheta$ and $g\in G$ with $g'(0)^q=1$
for some integer $q\geq1$. There is $K>0$, depending on $g$ and $S$, with
the following property: if $\rho>0$ is so small that $M=\nu(S,\rho)$
satisfies $1\leq M<\infty$ and $M\rho<K$, then $g$ has finite order.
\end{lemma}

\begin{proof}
Let $r_1$, $C$ and $\theta_c$ be as in the proof of
Lemma~\ref{lem:elliptic}, and $\lambda=g'(0)$, so that
$q\arg\lambda\in2\pi\Z$. Put $k=qM$ and
\[
W=\Bigl\{z:\ e\rho<|z|<\frac{\vartheta}{64\,CqM},\
|\arg z-\theta_c|<\frac\vartheta8\Bigr\},
\]
which is open and non-empty as soon as $M\rho<K=\vartheta/(64\,e\,Cq)$. For
$z\in W$ we have $|z|\leq1/(8kC)$, and the iteration estimate in the proof
of Lemma~\ref{lem:elliptic}, applied with $i=qj\leq k$, shows that the
points $z_j=g^{qj}(z)$, $0\leq j\leq M$, belong to the class of $z$, have
modulus in $(\rho,r)$, and have arguments which differ from $\theta_c$ by
at most $\vartheta/8+8Ck|z|\leq\vartheta/4$, since $qj\arg\lambda\in2\pi\Z$.
As in the proof of Lemma~\ref{lem:elliptic}, two of these $M+1$ points
coincide, $W$ is covered by the fixed-point sets of
$g^q,\dots,g^{qM}$, and some $g^{qd}$, $1\leq d\leq M$, is the identity.
\end{proof}

\begin{theorem}\label{thm:trichotomy}
Let $\sim$ be an equivalence relation on a neighbourhood of $0\in\C$, let
$S$ be a sector with vertex $0$, and let $g\in\Inv(\sim)$ be an element of
infinite order, with multiplier $\lambda=g'(0)$. There are $c>0$ and
$\rho_0>0$ such that, for $0<\rho<\rho_0$:
\begin{enumerate}[label=\textup{(\alph*)}]
\item if $|\lambda|\neq1$, then $\nu(S,\rho)\geq c\log(1/\rho)$;
\item if $\lambda$ is a root of unity, then $\nu(S,\rho)\geq c\,\rho^{-1}$;
\item if $|\lambda|=1$ and $\lambda$ is not a root of unity, then
$\nu(S,\rho)\geq c\,\rho^{-1/2}$.
\end{enumerate}
The exponents in \textup{(a)} and \textup{(b)} cannot be improved.
\end{theorem}

\begin{proof}
For $\rho$ smaller than the radius of $S$, $\nu(S,\rho)\geq1$. Part (a) is
Lemma~\ref{lem:hyperbolic}. In (b), Lemma~\ref{lem:parabolic} shows that,
for every small $\rho$, either $\nu(S,\rho)=\infty$ or
$\nu(S,\rho)\geq K/\rho$. In (c), Lemma~\ref{lem:elliptic} shows that
either $\nu(S,\rho)=\infty$ or $\nu(S,\rho)^2\geq K/\rho$.

For the optimality of (a), the relation of Theorem~\ref{thm:onevar}(c) has
$\nu(S,\rho)\leq\log_2(1/\rho)+1$ for every sector, and contains in its
invariance group the germ $z\mapsto z/2$. For (b), consider on the unit
disc the relation whose classes are $\{0\}$ and the sets of points $w\neq0$
with $1/w-1/z\in\Z$. Its invariance group contains the parabolic germ
$h(z)=z/(1+z)$, since $1/h(z)=1/z+1$, which has multiplier $1$ and
infinite order. The points of the class of $z$ with $|w|>\rho$ correspond
to the integers $n$ with $|1/z+n|<1/\rho$, so that
$\nu(S,\rho)\leq2/\rho+1$ for every sector $S$.
\end{proof}

The exponent in (c) comes from Dirichlet's theorem and is probably not
optimal. When $g$ is linearizable, its orbits are dense in circles, and
$\nu(S,\rho)=\infty$ for small $\rho$; the remaining case is that of the
germs of Cremer type, whose orbits are governed by the hedgehogs of
P\'erez-Marco \cite{PM97}.

\paragraph{Finitely generated groups.}
Let $G\subset\Diff(\C,0)$ be generated by $g_1,\dots,g_k$. A
\emph{realization} of $G$ consists of a disc $D=\{|z|<r_0\}$ and of
injective holomorphic maps $\hat g_i$ on $D$, fixing $0$, whose germs are
the $g_i$; the compositions of the maps $\hat g_i^{\pm1}$ whose successive
images lie in $D$ form a pseudogroup $\Gamma$, whose orbits define a
relation $\sim_\Gamma$ on $D$, with counts $\nu_\Gamma(S,\rho)$. These
counts measure how many points of an orbit stay at a given distance from
the fixed point; they should not be confused with the word growth of $G$.
A word in the $g_i^{\pm1}$ is defined, with all its successive images in
$D$, on a small disc about $0$, where it coincides with the corresponding
element of $G$; hence $G\subset\Inv(\sim_\Gamma)$, and
Theorem~\ref{thm:trichotomy}, together with Lemma~\ref{lem:finite},
gives the following.

\begin{corollary}\label{cor:fg}
Let $G\subset\Diff(\C,0)$ be a finitely generated group, $\Gamma$ the
pseudogroup of a realization of $G$, and $S$ a sector with vertex $0$.
Then $G$ is a finite cyclic group, analytically conjugate to a group of
rotations, in each of the following cases:
\begin{enumerate}[label=\textup{(\roman*)}]
\item $\liminf_{\rho\to0}\nu_\Gamma(S,\rho)/\log(1/\rho)=0$;
\item the multipliers of the generators have modulus one, and
$\liminf_{\rho\to0}\rho^{1/2}\,\nu_\Gamma(S,\rho)=0$;
\item the multipliers of the generators are roots of unity, and
$\liminf_{\rho\to0}\rho\,\nu_\Gamma(S,\rho)=0$.
\end{enumerate}
The thresholds in \textup{(i)} and \textup{(iii)} are attained by the
groups generated by $z/2$ and by $z/(1+z)$.
\end{corollary}

\begin{proof}
Since the multiplier map is a homomorphism, in case (ii), respectively
(iii), all the multipliers of $G$ have modulus one, respectively are
roots of unity. By Theorem~\ref{thm:trichotomy}, applied to
$\Inv(\sim_\Gamma)\supset G$, every element of $G$ has finite order; by
Lemma~\ref{lem:finite}, the multiplier map is then injective on $G$,
whose image is a finitely generated group of roots of unity, and a finite
subgroup of $\Diff(\C,0)$ is conjugate to a group of rotations. The last
assertion follows from the computations in the proof of
Theorem~\ref{thm:trichotomy}.
\end{proof}

\subsection{Virtual holonomy}

Applied to the relation $\sim_U$ on a transverse section, whose
invariance group is the virtual holonomy group,
Theorem~\ref{thm:onevar}, Step~2 of its proof and
Theorem~\ref{thm:trichotomy} read as follows; recall
that the holonomy maps of the leaf through the base point belong to this
group (Lemma~\ref{lem:vir}(c)).

\begin{corollary}\label{cor:leafAB}
Let $\Fcal_U$ be a representative of $\Fcal$, $\Sigma$ a transverse section
at $a$, and $S\subset\Sigma$ a sector of opening $\vartheta$.
\begin{enumerate}[label=\textup{(\alph*)}]
\item If $\nu_U(S)\leq N$, then $\Holv(\Fcal_U,\Sigma,a)$ is a finite
cyclic group of order at most $2\pi N/\vartheta$.
\item If $S$ has sub-logarithmic counts, then every element of
$\Holv(\Fcal_U,\Sigma,a)$ has finite order, the group is abelian, and
each of its finitely generated subgroups is finite and cyclic.
\item If $\Holv(\Fcal_U,\Sigma,a)$ contains an element of order $m$, then
$\nu_U(S)\geq\lceil m\vartheta/2\pi\rceil$.
\item If $\Holv(\Fcal_U,\Sigma,a)$ contains an element of infinite order
with multiplier $\lambda$, then there are $c>0$ and $\rho_0>0$ such that,
for $0<\rho<\rho_0$, $\nu_U(S,\rho)\geq c\log(1/\rho)$ if
$|\lambda|\neq1$, $\nu_U(S,\rho)\geq c\rho^{-1}$ if $\lambda$ is a root
of unity, and $\nu_U(S,\rho)\geq c\rho^{-1/2}$ otherwise.
\end{enumerate}
The same statements hold for a holomorphic foliation with isolated
singularities on a complex surface.
\end{corollary}

\section{Leaf counts and first integrals: proof of Theorem~\ref{thm:plane}}
\label{sec:foliations}

\subsection{Transport of sub-logarithmic counts}

In this subsection and in the next one, $\Gcal$ is a holomorphic foliation
with isolated singularities on a complex surface $X$. Leaves, transverse
sections, holonomy and virtual holonomy groups are taken with respect to
the leaves of $\Gcal$ on $X$ (Section~\ref{sec:virtual}), and we write
$\Holv(\Sigma,a)$ for $\Holv(\Gcal,\Sigma,a)$.

Transport by holonomy and through the corners of the reduction of
singularities was used in \cite[Lemmas~3.5--3.6 and
Proposition~3.9]{RS21} for moderate sectorial first integrals. The
following lemma transports sub-logarithmic leaf counts. At a singularity in
the linear normal form of~\ref{R3}, the passage from one separatrix to the
other is the Dulac correspondence \cite[\S2, pp.~372--373]{Sca02}; see
also \cite{CaSc00} and \cite[\S3.5]{RS21}.

\begin{lemma}[Transport of sub-logarithmic counts]\label{lem:transport}
\begin{enumerate}[label=\textup{(\alph*)}]
\item Let $C^\circ\subset X\setminus\sing(\Gcal)$ be a connected smooth
$\Gcal$-invariant curve. If some transverse section at a point of
$C^\circ$ contains a sector with sub-logarithmic counts, then every
transverse section $\Sigma$ at any point $x\in C^\circ$, with any
coordinate, contains such a sector, and every element of $\Holv(\Sigma,x)$
has finite order. The same holds with bounded counts.
\item Let $s$ be a non-degenerate reduced singularity of $\Gcal$, with
separatrices $C$ and $C'$, and suppose that some transverse section at a
point of $C\setminus\{s\}$ close to $s$ contains a sector with
sub-logarithmic counts. Then there are coordinates centred at $s$ in which
$\Gcal$ is defined by $p\,y\,dx+q\,x\,dy$, with $p,q\geq1$ coprime,
$C=\{y=0\}$ and $C'=\{x=0\}$, and every transverse section at a point of
$C'\setminus\{s\}$ close to $s$ contains a sector with sub-logarithmic
counts. The same holds with bounded counts.
\end{enumerate}
\end{lemma}

\begin{proof}
(a) The curve $C^\circ$ is an open subset of a leaf. Transport along a
path in $C^\circ$ gives a sector with sub-logarithmic, or bounded, counts,
possibly of smaller opening, by Lemma~\ref{lem:counts}(d), and
Corollary~\ref{cor:leafAB}(b) shows that the elements of the virtual
holonomy group have finite order.

(b) By (a), applied to a punctured disc of $C$ about $s$, the holonomy map
of a positive loop around $s$ in $C$ belongs to a virtual holonomy group
all of whose elements have finite order (Lemma~\ref{lem:vir}(c)); hence it
is periodic, and \ref{R3} gives the linear normal form, with $C=\{y=0\}$
and $C'=\{x=0\}$. The first integral $P=x^py^q$ is not a power, since $p$
and $q$ are coprime. Let $V$ be a neighbourhood of $s$, contained in the
domain of these coordinates, as in Lemma~\ref{lem:connected} for $g=P$.
Choose points $(x_0,0)$ and $(0,y_0)$ of $V$ with $x_0y_0\neq0$,
transverse sections $\Sigma_1\subset\{x=x_0\}\cap V$ and
$\Sigma_2\subset\{y=y_0\}\cap V$ at these points, and on them the
coordinates $v=x_0^{p/q}y$ and $w=y_0^{q/p}x$, for fixed choices of the
roots, so that $P|_{\Sigma_1}=v^q$ and $P|_{\Sigma_2}=w^p$. By (a) there
is a sector $S$ in $\Sigma_1$, in the coordinate $v$, with sub-logarithmic
counts. Choose a nonempty open interval $I$ of length at most $2\pi$ such
that $(p/q)I$ is compactly contained in the interval of arguments of $S$.
On a sufficiently small sector $S'$ of arguments $I$ in $\Sigma_2$, in the
coordinate $w$, the branch $\beta(w)=\exp((p/q)\log w)$ of the Dulac
correspondence takes its values in $S$. It is injective, because the
interval $(p/q)I$ has length less than $2\pi$, and $|v(\beta(w))|=|w|^{p/q}$.
Moreover $P(\beta(w))=\beta(w)^q=w^p=P(w)\neq0$, so $\beta(w)$ and $w$ are
on the same leaf of $\Gcal$, by Lemma~\ref{lem:connected}(a).
Lemma~\ref{lem:counts}(c), with $c=1$ and $\kappa=p/q$, shows that $S'$
has sub-logarithmic counts, or bounded counts if $S$ has, and (a) carries
this property to the other transverse sections to $C'\setminus\{s\}$.
\end{proof}

\subsection{The planar construction}

By an \emph{irreducible germ of invariant curve} at a point $s\in X$ we
mean an irreducible germ $C$ of analytic curve at $s$ which is
$\Gcal$-invariant. We denote by $C\setminus\{s\}$ the punctured part of a
representative of $C$ which is so small that $C\setminus\{s\}$ is a
connected smooth curve contained in $X\setminus\sing(\Gcal)$. By
Lemma~\ref{lem:transport}(a), the existence of a sector with
sub-logarithmic counts in a transverse section at a point of
$C\setminus\{s\}$ does not depend on the point, on the transverse section,
on its coordinate or on the representative.

\begin{proposition}\label{prop:planar}
Let $s$ be a singular point of\/ $\Gcal$ at which the germ of\/ $\Gcal$ is
a non-dicritical generalized curve, and let $C$ be an irreducible germ of
invariant curve at $s$. Suppose that some transverse section at a point of
$C\setminus\{s\}$ contains a sector with sub-logarithmic counts. Then:
\begin{enumerate}[label=\textup{(\alph*)}]
\item for every irreducible germ of invariant curve $C'$ at $s$, every
transverse section at a point of $C'\setminus\{s\}$ contains a sector with
sub-logarithmic counts;
\item the germ of\/ $\Gcal$ at $s$ has a non-constant holomorphic first
integral.
\end{enumerate}
\end{proposition}

The proof follows the proof of Th\'eor\`eme~1 of Mattei and Moussu
\cite[Chapitre~V, \S3]{MM80}, as in \cite[\S4.1]{RS21}: an induction on the
number of blow-ups, in which the local first integrals at the singular
points of the first blow-up are glued by means of
\cite[Chapitre~V, \S1, Lemme~3]{MM80}. The gluing lemma requires the
group generated by finitely many finite groups of germs to be finite. In
\cite{MM80} this follows from the topological hypothesis; here it follows
from Corollary~\ref{cor:leafAB}(b), since a finitely generated subgroup of
a virtual holonomy group all of whose elements have finite order is
finite.

\begin{proof}
We argue by induction on the number $r$ of blow-ups of the minimal
reduction of singularities of the germ of $\Gcal$ at $s$.

\emph{The case $r=0$.} The point $s$ is reduced, and it is non-degenerate
because the germ is a generalized curve. By~\ref{R1}, the irreducible germs
of invariant curves at $s$ are $C$ and the other separatrix $C'$. Hence
Lemma~\ref{lem:transport} gives (a), and (b) holds with the first integral
$x^py^q$, in the coordinates of Lemma~\ref{lem:transport}(b).

\emph{The case $r\geq1$: the first blow-up.} Let
$\pi_1\colon X_1\to X$ be the blow-up of $s$, with exceptional line
$E=\pi_1^{-1}(s)$, and let $\Gcal_1$ be the strict transform of $\Gcal$.
Since the germ of $\Gcal$ at $s$ is non-dicritical, $E$ is
$\Gcal_1$-invariant. Let $t_1,\dots,t_k$ be the singular points of
$\Gcal_1$ on $E$, and $E^\circ=E\setminus\{t_1,\dots,t_k\}$. The minimal
reduction of the germ of $\Gcal$ at $s$ is $\pi_1$ followed by the minimal
reductions of the germs of $\Gcal_1$ at the points $t_j$. Hence each of
these germs is a non-dicritical generalized curve whose minimal reduction
has fewer than $r$ blow-ups. The map $\pi_1$ restricts to a biholomorphism
$X_1\setminus E\to X\setminus\{s\}$ conjugating $\Gcal_1$ to $\Gcal$, and
$E^\circ$ is a leaf of $\Gcal_1$. Hence, for transverse sections contained
in $X_1\setminus E$, the leaves, the leaf counts and the virtual holonomy
groups of $\Gcal_1$ on $X_1$ correspond under $\pi_1$ to those of $\Gcal$
on $X$; we identify such a transverse section with its image.

Let $C'$ be an irreducible germ of invariant curve at $s$. Its strict
transform $\widetilde C'$ meets $E$ at a single point, as one sees by
lifting a Puiseux parametrization, and this point is one of the $t_j$,
because at a regular point of $\Gcal_1$ on $E$ the only germ of invariant
curve is the germ of $E$. Thus $\widetilde C'$ is an irreducible germ of
invariant curve of $\Gcal_1$ at $t_j$, different from $E$.

\emph{Transport along $E$ and proof of \textup{(a)}.} Let $t_1$ be the
point of $\widetilde C$. By the induction hypothesis~(a) at $t_1$, applied
to the invariant curve $\widetilde C$, a transverse section at a point of
$E^\circ$ contains a sector with sub-logarithmic counts. By
Lemma~\ref{lem:transport}(a), every transverse section $\tau$ at every
point $x\in E^\circ$ contains such a sector, and every element of
$\Holv(\Gcal_1,\tau,x)$ has finite order. Thus the induction hypothesis
applies at each point $t_j$, with the invariant curve $E$. Its part~(a),
at the point $t_j$ of $\widetilde C'$, gives (a) for $C'$.

\emph{Proof of \textup{(b)}.} By the induction hypothesis~(b), the germ of
$\Gcal_1$ at $t_j$ has a non-constant holomorphic first integral $g_j$,
$j=1,\dots,k$; by Lemma~\ref{lem:fi}(d) we may assume that $g_j(t_j)=0$
and that $g_j$ is not a power. We now use the notation of
\cite[Chapitre~V, \S1]{MM80}. Choose coordinates $(x,y)$ centred at $s$
such that the points of $E$ defined by the directions $\{y=0\}$ and
$\{x=0\}$ are not singular points of $\Gcal_1$, let $t_0\in E^\circ$ be the
first of these points, and consider the projection $(x,t)\mapsto t$ in the
chart $y=tx$ of the blow-up; its fibres over the points of $E^\circ$ are
transverse sections, and we denote by $\Sigma_0$ the fibre over $t_0$. Let
$V_j$ be a neighbourhood of $t_j$ as in Lemma~\ref{lem:connected} for
$g_j$, so small that $g_j$ vanishes on $E\cap V_j$
(Lemma~\ref{lem:fi}(a)), and let $\tau_j\subset V_j$ be the fibre over a
point $x_j\in E^\circ\cap V_j$. By Lemma~\ref{lem:connected}(b),
\[
\Inv(g_j|_{\tau_j})\subset\Holv(\Gcal_1,\tau_j,x_j).
\]
The germ $g_j^0$ of \cite{MM80} is obtained by transporting $g_j|_{\tau_j}$
to $\Sigma_0$ by holonomy along a path $c$ in $E^\circ$ from $x_j$ to
$t_0$. Hence its invariance group
$H(g_j^0)=\hol_c\circ\Inv(g_j|_{\tau_j})\circ\hol_c^{-1}$, which is finite
and cyclic by Lemma~\ref{lem:finite}(b), is contained in
$\Holv(\Gcal_1,\Sigma_0,t_0)$, by Lemma~\ref{lem:vir}(d). All the elements
of this group have finite order, so the group $H(g_1,\dots,g_k;t_0)$
generated by the finite groups $H(g_j^0)$, being finitely generated, is
finite (Corollary~\ref{cor:leafAB}(b)), and
\cite[Chapitre~V, \S1, Lemme~3, p.~502]{MM80} gives a non-constant
holomorphic first integral of the germ of $\Gcal$ at $s$.
\end{proof}

\subsection{Germs with a first integral}

Recall that a germ with a non-constant holomorphic first integral has a
primitive one (Lemma~\ref{lem:primitive}). The following proposition
contains Theorem~\ref{thm:plane}(c).

\begin{proposition}\label{prop:integrable-counts}
Let $\Fcal$ be a germ of singular holomorphic foliation at $0\in\C^2$
with a non-constant holomorphic first integral.
Let $F_0$ be a primitive first integral, $\Gamma$ a separatrix and
$m=\ord_\Gamma F_0$. Then:
\begin{enumerate}[label=\textup{(\alph*)}]
\item $N_\Gamma(\vartheta)=\lceil m\vartheta/2\pi\rceil$ for every
$\vartheta\in(0,2\pi]$;
\item if $\Fcal_U$ is a representative adapted to $\Gamma$ such that $F_0$
is defined on $U$, and $\Sigma$ is a transverse section at
$a\in\Gamma^\circ_U$, then $\Holv(\Fcal_U,\Sigma,a)=\Inv(F_0|_\Sigma)$, a
cyclic group of order $m$, and, in a coordinate $w$ such that
$F_0|_\Sigma=F_0(0)+w^m$, every sector $S\subset\Sigma$ of opening
$\vartheta$ and sufficiently small radius satisfies
$\nu_U(S)=\lceil m\vartheta/2\pi\rceil$.
\end{enumerate}
\end{proposition}

\begin{proof}
\emph{Upper bound.} Let $\Fcal_U$ be a representative adapted to $\Gamma$
such that $F_0$ is defined on $U$, let $\Sigma$ be a transverse section at
$a\in\Gamma^\circ_U$, and let $w$ be a coordinate on $\Sigma$ such that
$F_0|_\Sigma=F_0(0)+w^m$ on a disc $\Sigma'=\{|w|<\epsilon\}$
(Lemma~\ref{lem:fi}(b)). By Lemma~\ref{lem:fi}(a), points of $\Sigma'$ on
the same leaf of $\Fcal_U$ have the same value of $w^m$. Hence
$\Holv(\Fcal_U,\Sigma,a)\subset\Inv(F_0|_\Sigma)$, which is cyclic of order
$m$ by Lemma~\ref{lem:finite}(b), and the trace of a leaf on a sector
$S\subset\Sigma'$ is contained in a set $\{w^m=t\}$, $t\neq0$, which
consists of $m$ points equally spaced on a circle. By
Lemma~\ref{lem:arc}(a), if $S$ has opening $\vartheta$ in the coordinate
$w$, it contains at most $\lceil m\vartheta/2\pi\rceil$ of them. Hence
$\nu_U(S)\leq\lceil m\vartheta/2\pi\rceil$, and
$N_\Gamma(\vartheta)\leq\lceil m\vartheta/2\pi\rceil$.

\emph{Lower bound.} Let $\Fcal_U$ be any representative adapted
to $\Gamma$, $\Sigma$ a transverse section at $a\in\Gamma^\circ_U$, with
any coordinate, and $S\subset\Sigma$ a sector of opening $\vartheta$. By
Lemma~\ref{lem:primitive}, $F_0-F_0(0)$ is not a power. Let $\Fcal_{U'}$,
$U'\subset U$, be a representative adapted to $\Gamma$ such that $F_0$ is
defined on $U'$, let $V\subset U'$ be a neighbourhood of $0$ as in
Lemma~\ref{lem:connected}, for the foliation $\Fcal_U$ and the first
integral $g=F_0-F_0(0)$, which vanishes on $\Gamma$, and let $\Sigma'$ be a
transverse section at a point $a'\in\Gamma^\circ_U\cap V$; note that
$a'\in\Gamma^\circ_{U'}$. By Lemmas~\ref{lem:connected}(b)
and~\ref{lem:finite}(b), $\Holv(\Fcal_U,\Sigma',a')$ contains the group
$\Inv(F_0|_{\Sigma'})$, which is cyclic of order $m$
(Lemma~\ref{lem:fi}(b)). Since $\Gamma^\circ_U$ is a leaf of $\Fcal_U$
containing $a$ and $a'$, Lemma~\ref{lem:vir}(d) shows that
$\Holv(\Fcal_U,\Sigma,a)$ is conjugate to $\Holv(\Fcal_U,\Sigma',a')$.
Hence it contains an element of order $m$, and
$\nu_U(S)\geq\lceil m\vartheta/2\pi\rceil$ by
Corollary~\ref{cor:leafAB}(c). Together with the upper bound, this proves
(a) and (b): under the hypotheses of (b), the group
$\Holv(\Fcal_U,\Sigma,a)$ is contained in $\Inv(F_0|_\Sigma)$ and contains
an element of order $m$.
\end{proof}

\subsection{Saddle-nodes}

\begin{proposition}\label{prop:exposed}
Let $\pi\colon(M,D)\to(\C^2,0)$ be a reduction of singularities of a germ
of singular foliation $\Fcal$ at $0\in\C^2$, where the exceptional divisor
$D$ is empty if $\pi$ is the identity. Suppose that a saddle-node $s$ of
$\tF$ over $0$ has its strong separatrix $C_s$ not contained in $D$. Then
$\Gamma_s=\pi(C_s)$ is a separatrix of $\Fcal$, and no sector transverse
to $\Gamma_s$, in any representative adapted to $\Gamma_s$, has
sub-logarithmic counts. In particular $N_{\Gamma_s}\equiv\infty$.
\end{proposition}

\begin{proof}
Let $c$ be an injective holomorphic parametrization of a small
representative of the smooth curve $C_s$, with $c(0)=s$. Since
$C_s\not\subset D$, we may take this representative so small that
$c(t)\notin D$ for $t\neq0$. The map $\pi$ is a biholomorphism off
$\pi^{-1}(0)$ onto a punctured neighbourhood of $0$ and conjugates $\tF$
to $\Fcal$ there. Hence $\Gamma_s=\pi(C_s)$, the image of the non-constant
holomorphic germ $\pi\circ c$, is an irreducible analytic curve germ at
$0$ \cite{Chi89}, and it is invariant at its points other than $0$; that
is, $\Gamma_s$ is a separatrix of $\Fcal$.

Let $\Fcal_U$ be adapted to $\Gamma_s$ and take a small positive loop on
$C_s\setminus\{s\}$ in $\pi^{-1}(U)$. Since $\pi$ is a biholomorphism off
$\pi^{-1}(0)$, it takes a chain of flow boxes along this loop to a chain
along its image, so the two holonomy maps are conjugate by the restriction
of $\pi$ to the transverse section. The image of the loop is contained in
$\Gamma^\circ_{s,U}$, and its holonomy map has infinite order,
by~\ref{R2}. By Lemma~\ref{lem:vir}(c),(d), the virtual holonomy group of
$\Fcal_U$ at every transverse section at a point of $\Gamma^\circ_{s,U}$
contains an element of infinite order, and Corollary~\ref{cor:leafAB}(b)
shows that no sector of such a section has sub-logarithmic counts.
\end{proof}

\subsection{Proof of Theorem~\ref{thm:plane}}

\emph{Part \textup{(c)}} is Proposition~\ref{prop:integrable-counts}(a).

\emph{Part \textup{(b)}.} Let $\Fcal_U$ be a representative adapted to
$\Gamma$ and $S$ a sector transverse to $\Gamma$ with sub-logarithmic
counts. Proposition~\ref{prop:planar}(b), applied to $\Fcal_U$ on $X=U$,
with $s=0$ and $C=\Gamma$, gives the first integral.

\emph{Part \textup{(a)}.} If $\Fcal$ has a non-constant holomorphic first
integral, Proposition~\ref{prop:integrable-counts}(b) gives, for every
separatrix, sectors with bounded, hence sub-logarithmic, counts.
Conversely, assume that every separatrix $\Gamma$ of $\Fcal$ has a sector
with sub-logarithmic counts. We use the index theorem of Camacho and Sad
\cite[Section~2]{CS82}. Let $E$ be a smooth compact $\tF$-invariant curve.
For a singular point $s\in E$ of $\tF$, with local coordinates $(u,v)$ in
which $E=\{v=0\}$ and a generator $A(u,v)\partial_u+B(u,v)\partial_v$ of
$\tF$, the index of $\tF$ along $E$ at $s$ is
$\CS_s(\tF,E)=\Res_{u=0}\bigl(\partial_vB(u,0)/A(u,0)\bigr)\,du$, and the
sum of these indices over the singular points $s\in E$ equals the
self-intersection number $E\cdot E$. At a non-degenerate singularity the
index is the quotient of the eigenvalue transverse to $E$ by the
eigenvalue tangent to $E$; in particular, at a non-degenerate corner
$E_i\cap E_j$ the two indices are nonzero and inverse to each other. Along
the strong separatrix of a saddle-node the index is zero, since the
numerator vanishes at the origin and the denominator has a simple zero.

Let $\pi$ be the minimal reduction of singularities of $\Fcal$, with
strict transform $\tF$ and exceptional divisor $D$. A saddle-node whose
strong separatrix is not contained in $D$ is impossible: by
Proposition~\ref{prop:exposed} it would give a separatrix $\Gamma_s$ of
$\Fcal$ without sectors with sub-logarithmic counts. If $D$ is empty, the
same argument excludes a saddle-node at the origin. In the rest of the
argument we suppose that $D$ is nonempty and that $\tF$ has a saddle-node,
and we derive a contradiction.

\emph{A connected union of components bounded by strong separatrices.} The
choice of components and the intersection-matrix argument below are
adapted from \cite[pp.~314--316]{Tom99}; see also \cite[Theorem~3.1 and
Corollary~3.8]{San20}. Every component of $D$ is invariant, since $\Fcal$
is non-dicritical. Delete from the dual tree of $D$ the edges
corresponding to the corners which are saddle-nodes, and contract each
connected subtree of the resulting forest to a vertex; the deleted edges
are the edges of a new tree. Orient each such edge from the vertex
containing the component which is the strong separatrix of the
saddle-node to the vertex containing the component which is its weak
separatrix; these components are in distinct subtrees because the original
graph is a tree. Choose a vertex with no incoming edge and let $B$ be the
union of the corresponding components of $D$. Each corner between $B$ and
$D\setminus B$ is a saddle-node whose strong separatrix lies in $B$. Any
saddle-node lying on $B$ which is not a corner of $D$ also has its strong
separatrix in $B$, because the strong separatrices not contained in $D$
were already excluded. Moreover, at least one saddle-node has its strong
separatrix in $B$: if the new tree has more than one vertex, the chosen
vertex has an outgoing edge; if it has one vertex, then $B=D$ and there is
a saddle-node by assumption.

\emph{A separatrix of $\Fcal$ meeting $B$ at a non-degenerate point.}
Suppose that no non-degenerate singularity of $\tF$ on $B$ lies outside
the corners of $D$. For each component $E_i$ of $B$, the index theorem
then gives
\[
 E_i\cdot E_i=\sum_{j\sim i} c_{ij},\qquad
 c_{ij}=\CS_{E_i\cap E_j}(\tF,E_i),\qquad c_{ij}c_{ji}=1 ,
\]
where $j\sim i$ ranges over the components of $B$ meeting $E_i$. Indeed,
the saddle-nodes on $B$ contribute zero along their strong separatrices,
and every corner between two components of $B$ is non-degenerate. Choose a
root component and assign to it the weight $a_i=1$. Recursively assign
$a_j=-c_{ij}a_i$ across each edge directed away from the root of this
subtree. All the weights are nonzero, and the relations
$c_{ij}c_{ji}=1$ give $a_j=-c_{ij}a_i$ in both directions. Thus
\[
 (E_i\cdot E_i)\,a_i+\sum_{j\sim i}a_j=0
\]
for every component $E_i$ of $B$. This places the nonzero vector $(a_i)$
in the kernel of the intersection matrix of the components of $B$. But the
intersection matrix of $D$ is negative definite: in the basis of the total
transforms of the exceptional curves of the successive blow-ups the matrix
is diagonal with entries $-1$, and the strict transforms form another
basis. The matrix of $B$ is a principal submatrix, hence invertible, a
contradiction. Consequently some non-degenerate singular point of $\tF$ on
$B$ is not a corner of $D$. Its second separatrix is smooth and transverse
to $D$ by~\ref{R1}, and its image under $\pi$ is a separatrix $\Gamma$ of
$\Fcal$.

\emph{Transport to a strong separatrix.} By hypothesis there are a
representative $\Fcal_U$ adapted to $\Gamma$ and a sector transverse to
$\Gamma$ with sub-logarithmic counts. On $M_U=\pi^{-1}(U)$ the leaves of
$\tF$ off $D$ correspond to those of $\Fcal_U$ off $0$; in particular, the
lifted sector has sub-logarithmic counts for the leaves of $\tF$ on $M_U$.
Transport this property along the strict transform of $\Gamma$ to its
intersection with $B$, and through that non-degenerate singularity, by
Lemma~\ref{lem:transport}. Applying the same lemma along the components of
$B$ and through the corners between them, which are non-degenerate, we
obtain that every element of the virtual holonomy groups at the
transverse sections to the regular part of each component of $B$ has
finite order. This includes a transverse section to a component which is
the strong separatrix of a saddle-node. The holonomy map of a loop around
the saddle-node in that strong separatrix has infinite order, by~\ref{R2},
a contradiction.

Therefore the reduction of singularities has no saddle-nodes, and $\Fcal$
is a non-dicritical generalized curve. It has a separatrix by \cite{CS82},
with a sector with sub-logarithmic counts by assumption, and part~(b)
gives a non-constant holomorphic first integral. \qed

\subsection{Proof of the corollaries}

\begin{proof}[Proof of Corollary~\textup{\ref{cor:gap}}]
If $\Fcal$ has a non-constant holomorphic first integral, (i) holds by
Theorem~\ref{thm:plane}(c), and (ii) fails because the sectors of
Proposition~\ref{prop:integrable-counts}(b) have bounded counts. If
$\Fcal$ has no such first integral, (i) fails, and no sector transverse to
$\Gamma$ has sub-logarithmic counts, by Theorem~\ref{thm:plane}(b); that
is, (ii) holds.
\end{proof}

\begin{proof}[Proof of Corollary~\textup{\ref{cor:modA}}]
Let $\varphi$ be a moderate sectorial first integral of a representative
$\Fcal_U$ adapted to a separatrix $\Gamma$, on a sector $S$ transverse to
$\Gamma$. As observed in Section~\ref{sec:sectors}, $\varphi$ has a
leading term $a_dz^d$ with $d\geq1$. By Lemma~\ref{lem:leading},
$\varphi$ is injective on a smaller sector $S'\subset S$, so that every
leaf of $\Fcal_U$ meets $S'$ in at most one point; in particular $S'$
has bounded, hence sub-logarithmic, counts. The two cases of
Corollary~\ref{cor:modA} now follow from Theorem~\ref{thm:plane}(b)
and~(a).
\end{proof}

\subsection{Examples}

\begin{example}[Sharpness of the counting formula]\label{ex:sharp}
Let $m\geq1$ and let $\Fcal_m$ be defined at $0\in\C^2$ by
$m\,y\,dx+x\,dy$. The first integral $P=x^my$ is not a power, hence it is
primitive (Lemma~\ref{lem:primitive}); it has order $m$ along
$\Gamma=\{x=0\}$ and order $1$ along $\{y=0\}$. By
Theorem~\ref{thm:plane}(c),
$N_\Gamma(\vartheta)=\lceil m\vartheta/2\pi\rceil$ and
$N_{\{y=0\}}\equiv1$. Given $N\geq1$ and $\vartheta\in(0,2\pi]$, take
$m=\lfloor2\pi N/\vartheta\rfloor$. By
Proposition~\ref{prop:integrable-counts}(b), a transverse section to
$\Gamma$ then contains sectors $S$ of opening $\vartheta$ with
$\nu_U(S)=\lceil m\vartheta/2\pi\rceil\leq N$, while its virtual holonomy
group and the primitive first integral have order
$m=\lfloor2\pi N/\vartheta\rfloor$ along $\Gamma$. Thus the bounds of
Theorem~\ref{thm:onevar}(a) and Theorem~\ref{thm:plane}(c) are attained.
\end{example}

\begin{example}[The hypotheses of Theorem~\ref{thm:plane}(a)]
\label{ex:necessity}
(a) \emph{A dicritical germ.} Let $\Fcal$ be the radial foliation defined
by $\omega=x\,dy-y\,dx$, and $U$ a ball. The leaves of $\Fcal_U$ are the
punctured lines through $0$, and every line through $0$ is a separatrix.
The transverse section $\{x=x_0\}$ to $\Gamma=\{y=0\}$ meets each leaf at
most once, so $N_\Gamma\equiv1$. A holomorphic first integral is constant
on each line through $0$, hence equal to its value at $0$. The germ
$\phi(u,v)=(v,uv)$ satisfies $\phi^*\omega=v^2\,du$ and maps $\{v=0\}$ to
the origin, so $\Fcal$ is dicritical.

(b) \emph{A non-dicritical germ which is not a generalized curve.} Let
$\Fcal$ be the saddle-node $y\,dx-x^{k+1}dy$, $k\geq1$, with strong
separatrix $\{x=0\}$ and weak separatrix $\Gamma=\{y=0\}$. On $\{x\neq0\}$
the function $H=y\exp\bigl(1/(kx^k)\bigr)$ satisfies
$dH=-\exp\bigl(1/(kx^k)\bigr)x^{-k-1}\bigl(y\,dx-x^{k+1}dy\bigr)$, so it
is a first integral there. Let $U$ be a bidisc and $\Sigma=\{x=x_0\}$,
$x_0\neq0$. A leaf meeting $\Sigma$ lies in $\{x\neq0\}$, because
$\{x=0\}$ is invariant, hence in a level of $H$; since
$H|_\Sigma=y\exp(1/(kx_0^k))$ is injective, every leaf meets $\Sigma$ at
most once, and $N_\Gamma\equiv1$. However, $\Fcal$ has no non-constant
holomorphic first integral, by~\ref{R2}; compare
\cite[Example~3.4]{RS21}. Thus, in Theorem~\ref{thm:plane}(a), the
condition cannot be imposed on one separatrix only; by
Proposition~\ref{prop:exposed}, it fails for the strong separatrix.
\end{example}

\begin{example}[Finite traces with counts of order $1/\rho$]\label{ex:finite}
Let $\Fcal$ be generated by $X=x\partial_x-y(1-xy)\partial_y$. The origin
is a non-degenerate reduced singularity with eigenvalues $1,-1$, so
$\Fcal$ is a non-dicritical generalized curve. Since $X(xy)=x^2y^2$ and
$X(\log x)=1$, the function $H=x\exp\bigl(1/(xy)\bigr)$ is a first
integral on the complement of the axes. Let $U$ be a bidisc,
$\Sigma=\{x=x_0\}$ with $x_0\neq0$, and $s=1/(x_0y)$. The leaves meeting
$\Sigma\setminus\{y=0\}$ stay in the complement of the axes, which are
invariant, and $H|_\Sigma=x_0e^s$; hence each leaf meets $\Sigma$ in a
subset of a set $s_0+2\pi i\Z$. Fix $0<\beta<\pi/2$ and let
$S\subset\Sigma$ be the sector $|\arg(x_0y)|<\pi/2-\beta$, $0<|y|<r$. Its
image in the $s$-plane lies in the cone $|\Im s|<(\cot\beta)\Re s$, which
contains finitely many points of each set $s_0+2\pi i\Z$, so every leaf
meets $S$ in finitely many points; those of modulus $|y|>\rho$ correspond
to points with $|s|<1/(|x_0|\rho)$, so that $\nu_U(S,\rho)$ is at most of
order $1/\rho$. On the other hand, following a positive loop
$x=x_0e^{i\theta}$ along $\{y=0\}$ and using the constancy of
$\log x+1/(xy)$ on the leaves, one finds that the holonomy map of
$\{y=0\}$ is $s\mapsto s-2\pi i$, that is, $y\mapsto y/(1-2\pi ix_0y)$, a
parabolic germ of infinite order. By Lemma~\ref{lem:vir}(c),(d) and
Corollary~\ref{cor:leafAB}(b), in every representative adapted to
$\{y=0\}$ no sector transverse to $\{y=0\}$ has sub-logarithmic counts;
for the sector $S$ above, the points $s_0+2\pi ik$ with $|k|$ up to a
constant times $1/\rho$ show that $\nu_U(S,\rho)$ is exactly of order
$1/\rho$. Thus every trace $L\cap S$ is finite, although $\nu_U(S)=\infty$
and $\Fcal$ has no non-constant holomorphic first integral.

The same germ shows that, in Corollary~\ref{cor:modA}, the asymptotic
expansion must be non-constant, and not merely nonzero. On the sector
$S^-\subset\Sigma$ given by $|\arg(x_0y)-\pi|<\pi/2-\beta$, $0<|y|<r$, we
have $\Re s\leq-\sin\beta/(|x_0||y|)$, so that $\varphi=1+x_0e^s$ is a
sectorial first integral with the nonzero constant expansion
$\widehat\varphi=1$. The same holds on every transverse section at a point
$(x_1,0)$, $x_1\neq0$, of every representative contained in a bidisc,
since there $1/(xy)$ has a simple pole in a coordinate centred at
$(x_1,0)$. Hence the hypotheses of \cite[Theorem~A]{RS21}, with a nonzero
asymptotic expansion, are satisfied, although its conclusion fails.
\end{example}

\begin{example}[Logarithmic counts]\label{ex:log}
Let $\alpha\in\C\setminus\R$ and let $\Fcal$ be defined by
$y\,dx-\alpha x\,dy$, that is, by the vector field
$\alpha x\partial_x+y\partial_y$. Its eigenvalues have the non-real
quotient $\alpha$, so the origin is a non-degenerate reduced singularity
and $\Fcal$ is a non-dicritical generalized curve. Let $U$ be the bidisc
$\{|x|<1,|y|<1\}$, $\Gamma=\{y=0\}$, and $\Sigma=\{x=x_0\}$ with
$0<|x_0|<1$, with the coordinate $y$. The leaves off the axes are given by
$y=c\,x^{1/\alpha}$, and the holonomy of $\Gamma$ along the loop
$x=x_0e^{i\theta}$ is $y\mapsto\mu y$, $\mu=e^{2\pi i/\alpha}$, with
$|\mu|\neq1$. A leaf meets $\Sigma$ in a subset of an orbit
$\{\mu^ky_0\}_{k\in\Z}$, so that
$\nu_U(S,\rho)\leq\log(1/\rho)/|\log|\mu||+1$ for every sector $S$ of
$\Sigma$. Conversely, $\mu$ belongs to the virtual holonomy group, and
Lemma~\ref{lem:hyperbolic} gives $\nu_U(S,\rho)\geq c\log(1/\rho)$ for
small $\rho$. Thus the leaf counts grow exactly logarithmically, and
$\Fcal$ has no non-constant holomorphic first integral: otherwise, by
Proposition~\ref{prop:integrable-counts}(b), the virtual holonomy group
of $\Gamma$ in a smaller representative would be finite, whereas it
contains the holonomy map $y\mapsto\mu y$, of infinite order. This shows that the threshold in
Theorem~\ref{thm:plane} and Corollary~\ref{cor:gap} is attained.
\end{example}
\section{Closed leaves and periodic holonomy: proofs of Theorems~\ref{thm:Q52} and~\ref{thm:group}}
\label{sec:Q52}

Throughout this section $F=xy(y-x)$,
$\Delta=a_0x^2+a_1xy+a_2y^2$ is a homogeneous quadratic polynomial,
$\lambda(\Delta)=2a_0+a_1+2a_2$, and $\Fcal_\Delta$ is the germ at
$0\in\C^2$ of the foliation defined by
\[
\omega_\Delta=\tfrac13\,dF+\Delta\,(x\,dy-y\,dx).
\]
We put $\Gamma_1=\{y=0\}$, $\Gamma_2=\{y=x\}$, $\Gamma_3=\{x=0\}$,
$P(t)=t(t-1)$, $\zeta=e^{2\pi i/3}$, $\delta(t)=\Delta(1,t)$ and
$\hat\delta(\hat t)=\Delta(\hat t,1)$. The germ $\Fcal_\Delta$ is reduced
by one blow-up, and the holonomy group of the exceptional line, which is
the \emph{projective holonomy} group of $\Fcal_\Delta$ in the terminology
of \cite{Sca16} (the holonomy of $\omega_\Delta$ in
\cite[Chapitre~III, \S5]{MM80}), is a group of affine maps in a suitable
transverse coordinate. The form $\omega_\Delta$ may have zeros other than
the origin; for $\Delta=\lambda xy$, $\lambda\neq0$, it has exactly one,
the point $\bigl(-1/(3\lambda),-1/(3\lambda)\bigr)$ of $\Gamma_2$. Such
zeros are excluded by the hypothesis on $U$ in Theorem~\ref{thm:Q52}.

As explained in the introduction, the theorem of Mattei and Moussu
\cite[Th\'eor\`eme~B]{MM80} does not apply because the leaves accumulate
at the separatrices; as in Moussu's example, the obstruction comes from a
transversely Euclidean structure with discrete monodromy, in the spirit
of \cite{CLNLPT07}.

\subsection{The blow-up and the projective holonomy}\label{sec:blowup}

The initial homogeneous part $\frac13\,dF$ of $\omega_\Delta$ has
relatively prime coefficients; hence $\omega_\Delta$ is saturated, and its
zero at the origin is isolated. In the chart $y=tx$ of the blow-up $\pi$
of the origin we have $F=x^3P(t)$, $x\,dy-y\,dx=x^2\,dt$ and
$\Delta(x,tx)=x^2\delta(t)$, hence
\[
\pi^*\omega_\Delta=x^2\Bigl(P(t)\,dx
+x\bigl(\tfrac13\,P'(t)+x\,\delta(t)\bigr)\,dt\Bigr).
\]
The strict transform $\tF_\Delta$ is defined by the form in the large
parentheses. The exceptional line $E=\{x=0\}$ is invariant, and on its
complement the transverse coordinate $W=1/x$ turns the equation of the
leaves into the linear equation
\begin{equation}\label{eq:affine}
\frac{dW}{dt}=\frac{P'(t)}{3P(t)}\,W+\frac{\delta(t)}{P(t)},\qquad
\frac{P'(t)}{3P(t)}=\frac13\Bigl(\frac1t+\frac1{t-1}\Bigr).
\end{equation}
In the chart $x=\hat t\,y$ we have $F=-y^3P(\hat t)$,
$x\,dy-y\,dx=-y^2\,d\hat t$ and $\Delta(\hat ty,y)=y^2\hat\delta(\hat t)$,
hence
\[
\pi^*\omega_\Delta=-y^2\Bigl(P(\hat t)\,dy
+y\bigl(\tfrac13\,P'(\hat t)+y\,\hat\delta(\hat t)\bigr)\,d\hat t\Bigr),
\]
and $\widehat W=1/y$ satisfies equation~\eqref{eq:affine} in the variable
$\hat t$, with $\hat\delta$ in place of $\delta$. On the overlap,
$\hat t=1/t$ and $\widehat W=W/t$. Thus the singular points of
$\tF_\Delta$ on $E$ are the points $s_1$, $s_2$ and $s_3$ given by $t=0$,
$t=1$ and $\hat t=0$, and $s_i$ lies on the strict transform of
$\Gamma_i$.

\begin{lemma}\label{lem:saddles}
Let $t_1=0$ and $t_2=1$ in the chart $y=tx$, and $t_3=0$ in the chart
$x=\hat t\,y$.
\begin{enumerate}[label=\textup{(\alph*)}]
\item Near $t_i$, the equation~\eqref{eq:affine} of the corresponding
chart has a unique holomorphic solution $\gamma_i$, and its solutions are
the functions $\gamma_i+c\,\varrho$, $c\in\C$, where $\varrho$ is a
branch of $P^{1/3}$.
\item Near $s_i$, $i=1,2,3$, there are coordinates $(\xi,\tau)$ in which
$E=\{\xi=0\}$, the strict transform of\/ $\Gamma_i$ is $\{\tau=0\}$, and
$\tau\xi^3$ is a first integral of\/ $\tF_\Delta$.
\item Fix $t_*\in\C\setminus\{0,1\}$. The holonomy maps of the leaf
$E^\circ=E\setminus\{s_1,s_2,s_3\}$ of\/ $\tF_\Delta$, computed on the
transverse section $\{t=t_*\}$, are given in the coordinate $W$ by the
monodromy of equation~\eqref{eq:affine}; for the positive loops around
$t_1$ and $t_2$ they have the form $f_i(W)=\zeta W+\mu_i$, $i=1,2$.
\end{enumerate}
\end{lemma}

\begin{proof}
(a) Write $P(t)=(t-t_i)\pi_i(t)$, with $\pi_i(t_i)\neq0$. A branch
$\varrho$ of $P^{1/3}$ near $t_i$ is a solution of the homogeneous
equation, and a solution of~\eqref{eq:affine} is the product of $\varrho$
by a primitive of $\delta/(P\varrho)=(t-t_i)^{-4/3}g_i(t)$, where $g_i$ is
holomorphic at $t_i$ (with $\hat\delta$ in place of $\delta$ for $i=3$).
Termwise integration of the expansion of $g_i$ in powers of $t-t_i$
involves the denominators $j-\frac13$, $j\geq0$, which do not vanish; it
gives a primitive of the form $(t-t_i)^{-1/3}k_i(t)$ with $k_i$
holomorphic at $t_i$, and hence the holomorphic solution
$\gamma_i=\pi_i^{1/3}k_i$. The other solutions are $\gamma_i+c\,\varrho$,
$c\neq0$, and they are ramified at $t_i$.

(b) By (a), the function $(W-\gamma_i(t))/\varrho(t)$ is constant on the
leaves near $s_i$, $i=1,2$. Hence, with $\xi=x/(1-x\gamma_i(t))$, that
is, $1/\xi=W-\gamma_i(t)$, and $\tau=P(t)$, the function
$\tau\xi^3=\varrho^3/(W-\gamma_i)^3$ is a first integral, and $(\xi,\tau)$
are local coordinates at $s_i$ in which $E=\{\xi=0\}$ and the strict
transform of $\Gamma_i$ is $\{t=t_i\}=\{\tau=0\}$. At $s_3$ the same
construction in the chart $x=\hat t\,y$, with $\xi=y/(1-y\gamma_3(\hat
t))$ and $\tau=P(\hat t)$, gives the result.

(c) Off $E$, the leaves of $\tF_\Delta$ over $\C\setminus\{0,1\}$ are the
graphs of the solutions of~\eqref{eq:affine}, and these solutions can be
continued along every path in $\C\setminus\{0,1\}$. Hence, if $c$ is a
loop in $\C\setminus\{0,1\}$ based at $t_*$, the holonomy map of $E^\circ$
along $c$ is given by the continuation of the solutions along $c$, which
is an affine map $f_c$ of the coordinate $W$; in the coordinate $x=1/W$ it
is a germ of biholomorphism at $x=0$. By (a), the continuation along a
small positive loop around $t_i$ fixes $\gamma_i$ and multiplies
$\varrho$ by $\zeta$; hence the linear part of $f_i$ is $\zeta$.
\end{proof}

\begin{proof}[Proof of Theorem~\textup{\ref{thm:Q52}(a)}]
By Lemma~\ref{lem:saddles}(b), near each $s_i$ the foliation
$\tF_\Delta$ is defined by $3\tau\,d\xi+\xi\,d\tau$: its singular points
are linearizable saddles with eigenvalue ratio $-3$, hence non-degenerate
and reduced, and $E$ is invariant. As a check, the Camacho--Sad index of
$\tF_\Delta$ along $E$ at each $s_i$ is $-1/3$, and the three indices
add up to $E\cdot E=-1$. Thus $\Fcal_\Delta$ is a non-dicritical
generalized curve which is reduced by one blow-up. The strict transform of
a separatrix of $\Fcal_\Delta$ meets $E$ at a singular point, where it is
the separatrix other than $E$, as in the proof of
Proposition~\ref{prop:planar}; hence the separatrices are $\Gamma_1$,
$\Gamma_2$ and $\Gamma_3$. Finally, on the transverse section
$\{\xi=\xi_0\}$ to $\{\tau=0\}$, with the coordinate $\tau$, the lift of
the loop $\xi=\xi_0e^{i\theta}$, $0\leq\theta\leq2\pi$, through the point
$\tau$ is $\tau e^{-3i\theta}$, because $\tau\xi^3$ is constant on the
leaves. Hence the holonomy map of this loop, which generates the
fundamental group of the punctured separatrix, is the identity.
\end{proof}

The criterion in Theorem~\ref{thm:Q52}(c) rests on the following lemma.

\begin{lemma}\label{lem:polsol}
The following are equivalent.
\begin{enumerate}[label=\textup{(\roman*)}]
\item $\lambda(\Delta)=0$.
\item There are $a,b\in\C$ with $3\Delta=a(x^2-2xy)+b(y^2-2xy)$.
\item Equation~\eqref{eq:affine} has a polynomial solution.
\item The maps $f_1$ and $f_2$ of Lemma~\textup{\ref{lem:saddles}(c)}
have a common fixed point.
\end{enumerate}
In this case $a$ and $b$ are unique, $W=a+bt$ is the polynomial solution,
and, with $\ell=ax+by$,
\begin{equation}\label{eq:exactfi}
\omega_\Delta=\tfrac13\,(1-\ell)^4\,d\Bigl(\frac{F}{(1-\ell)^3}\Bigr).
\end{equation}
\end{lemma}

\begin{proof}
(i)$\Leftrightarrow$(ii). The polynomials $x^2-2xy$ and $y^2-2xy$ are
linearly independent and satisfy $\lambda=0$, and $\lambda$ is a non-zero
linear form on the three-dimensional space of the $\Delta$; hence they span
its kernel, and $a=3a_0$, $b=3a_2$ are unique.

(ii)$\Leftrightarrow$(iii). Equation~\eqref{eq:affine} reads
$3PW'=P'W+3\delta$. For $W=bt+a$ it becomes
$3\delta(t)=bt^2-2(a+b)t+a$, that is,
$3\Delta(x,y)=3x^2\delta(y/x)=a(x^2-2xy)+b(y^2-2xy)$. A polynomial
solution has degree at most one, since a term of degree $d$ contributes
$(3d-2)\,t^{d+1}$ to $3PW'-P'W$, and $3\delta$ has degree at most two.

(iii)$\Rightarrow$(iv). A polynomial solution is single-valued, so its
value at $t_*$ is fixed by $f_1$ and $f_2$.

(iv)$\Rightarrow$(iii). The solution through a common fixed point is
single-valued on $\C\setminus\{0,1\}$. Near $t_i$, $i=1,2$, it then
coincides with the holomorphic solution $\gamma_i$, by
Lemma~\ref{lem:saddles}(a). For the same reason $\widehat W=W/t$, which is
a single-valued solution of the equation in the chart $x=\hat t\,y$, is
holomorphic at $\hat t=0$. Hence this solution is an entire function which
grows at most linearly at infinity, that is, a polynomial.

Finally, both sides of $\ell\,dF-3F\,d\ell=-3\Delta\,(x\,dy-y\,dx)$ vanish
on the radial vector field $x\partial_x+y\partial_y$, by Euler's identity,
and both take the value $3y\Delta$ on $\partial_x$, as a direct computation
using (ii) shows; hence the identity holds. Therefore
\[
d\Bigl(\frac{F}{(1-\ell)^3}\Bigr)
=\frac{(1-\ell)\,dF+3F\,d\ell}{(1-\ell)^4}
=\frac{dF-(\ell\,dF-3F\,d\ell)}{(1-\ell)^4}
=\frac{3\,\omega_\Delta}{(1-\ell)^4}. \qedhere
\]
\end{proof}

\begin{proof}[Proof of Theorem~\textup{\ref{thm:Q52}(c)}]
If $\lambda(\Delta)=0$, then $F/(1-ax-by)^3$ is a holomorphic first
integral, by~\eqref{eq:exactfi}, and it is not constant. Let
$\lambda(\Delta)\neq0$. By Lemma~\ref{lem:polsol}, the fixed points $c_1$
and $c_2$ of $f_1$ and $f_2$ are distinct, and
\[
f_1\circ f_2^{-1}(W)=W+(1-\zeta)(c_1-c_2)
\]
is a non-trivial translation. In the coordinate $x=1/W$ it is a holonomy
map $h$ of the leaf $E^\circ$ which is tangent to the identity and
different from the identity; hence $h$ has infinite order
(Lemma~\ref{lem:finite}).

Suppose that $\Fcal_\Delta$ has a non-constant holomorphic first integral
$H$. Then $H\circ\pi$ is a non-constant first integral of $\tF_\Delta$
on a neighbourhood of $E$. Its restriction $k$ to the transverse section
$\{t=t_*\}$ is not constant, since a first integral which is constant on a
transverse section is constant on a flow box, and $k\circ h=k$, because
$h$ moves the points within their leaves. This contradicts
Lemma~\ref{lem:finite}(b), according to which the invariance group of $k$
is finite.
\end{proof}

\begin{proof}[Proof of Theorem~\textup{\ref{thm:Q52}(d)}]
Let $\alpha\in\C^*$ and $p,q\in\C$, put
$\Psi(x,y)=(x,y)/(\alpha+px+qy)$, a germ of biholomorphism with
$D\Psi(0)=\alpha^{-1}\mathrm{Id}$, and
$\Delta'=\alpha\Delta+\frac13\bigl(p(x^2-2xy)+q(y^2-2xy)\bigr)$. The map
$\Psi$ preserves the lines through the origin; in the chart $y=tx$ of the
blow-up it preserves $t$ and acts on $W=1/x$ by
$W\mapsto\widetilde W=\alpha W+p+qt$. If $W$ is a solution
of~\eqref{eq:affine}, that is, $3PW'=P'W+3\delta$, then
\[
3P\widetilde W'-P'\widetilde W
=3\alpha\delta+3qP-P'(p+qt)
=3\alpha\delta+qt^2-2(p+q)t+p ,
\]
and the right-hand side is $3\delta'$, where $\delta'(t)=\Delta'(1,t)$.
Hence $\Psi$ maps the leaves of $\Fcal_\Delta$ contained in $\{x\neq0\}$
into leaves of $\Fcal_{\Delta'}$, so that
$\Psi^*\omega_{\Delta'}\wedge\omega_\Delta$ vanishes on $\{x\neq0\}$, hence
identically: $\Psi$ conjugates $\Fcal_\Delta$ to $\Fcal_{\Delta'}$.

By Lemma~\ref{lem:polsol}, the polynomials
$\frac13\bigl(p(x^2-2xy)+q(y^2-2xy)\bigr)$ form the kernel of $\lambda$.
If $\lambda(\Delta)=0$, choose $\alpha=1$ and $(p,q)$ with $\Delta'=0$.
If $\lambda(\Delta)\neq0$, choose $\alpha=1/\lambda(\Delta)$ and $(p,q)$
with $\Delta'=xy$; this is possible since $\lambda(xy)=1$. Finally, a
germ with a non-constant holomorphic first integral is not conjugate to
one without, so the two classes are distinct by~(c).
\end{proof}

\subsection{The projective holonomy group}\label{sec:group}

The holonomy group of the leaf $E^\circ$ of $\tF_\Delta$, computed on the
transverse section $\{t=t_*\}$, is the projective holonomy group of
$\Fcal_\Delta$. Since the fundamental group of $E^\circ$, a sphere minus
three points, is generated by positive loops around $t_1=0$ and $t_2=1$,
it is generated by the maps $f_1$ and $f_2$ of
Lemma~\ref{lem:saddles}(c). The following theorem describes it through a
one-variable model which depends on a single complex parameter.

\begin{theorem}\label{thm:group}
Let $\zeta=e^{2\pi i/3}$ and, for $\lambda\in\C$, let
$G_\lambda\subset\Diff(\C,0)$ be the group generated by
$g_1(x)=\zeta^2x$ and $g_{2,\lambda}(x)=x/(\zeta+\lambda x)$. For $r>0$,
consider the pseudogroup generated by $g_1$ and $g_{2,\lambda}$ on the
disc $D_r=\{|x|<r\}$.
\begin{enumerate}[label=\textup{(\alph*)}]
\item $g_1$ and $g_{2,\lambda}$ have order three for every $\lambda$.
\item For every $\lambda$ and every $r>0$, each orbit of the pseudogroup
on $D_r$ through a non-zero point is a closed subset of
$D_r\setminus\{0\}$, and there is a
non-constant $G_\lambda$-invariant meromorphic function on
$\C\setminus\{0\}$.
\item The following are equivalent: \textup{(i)} $\lambda=0$;
\textup{(ii)} $G_\lambda$ is finite; \textup{(iii)} there is a
non-constant $G_\lambda$-invariant germ of holomorphic function at $0$.
\item For $\lambda\neq0$, the homothety $x\mapsto\lambda x$ conjugates
$G_\lambda$ to $G_1$, which is isomorphic to the crystallographic group
$\Z^2\rtimes\Z/3\Z$. The group $G_\lambda$ contains the parabolic germ
$g_1\circ g_{2,\lambda}^{-1}(x)=x/(1-\lambda x)$.
\end{enumerate}
Moreover, for every homogeneous quadratic polynomial
$\Delta=a_0x^2+a_1xy+a_2y^2$, the projective holonomy group of the germ of
foliation at $0\in\C^2$ defined by
$\frac13\,d\bigl(xy(y-x)\bigr)+\Delta\,(x\,dy-y\,dx)$ is analytically
conjugate to $G_\lambda$, with $\lambda=2a_0+a_1+2a_2$.
\end{theorem}

Part~(b) is the case of order three in the list of discrete groups of
\cite[Proposition~4.21]{Lor21}; Moussu's cuspidal example \cite{Mou85}
realizes the case of order six \cite[Exemple~8.24]{Lor21}.

\begin{lemma}\label{lem:normal}
There are a constant $\kappa_0\neq0$ and, on $\{t=t_*\}$, an affine change
of coordinate $W'=(W-c_1(\Delta))/((1-\zeta)\kappa_0)$, where $c_1(\Delta)$
depends linearly on $\Delta$, in which $f_1(W')=\zeta W'$ and
$f_2(W')=\zeta W'+\lambda(\Delta)$. Consequently, in the coordinate
$x'=1/W'$ the projective holonomy group of $\Fcal_\Delta$ is the group
$G_{\lambda(\Delta)}$ of Theorem~\textup{\ref{thm:group}}.
\end{lemma}

\begin{proof}
By Lemma~\ref{lem:saddles}(c), $f_i(W)=\zeta W+\mu_i$, so $f_i$ is the
rotation of order three about its unique fixed point $c_i(\Delta)$. The
map $f_i$ is the monodromy of~\eqref{eq:affine} along a loop which follows
a path $\delta_i$ from $t_*$ to a point close to $t_i$, turns once around
$t_i$ and comes back along $\delta_i$; it fixes the value at $t_*$ of the
continuation along $\delta_i$ of the holomorphic solution $\gamma_i$ of
Lemma~\ref{lem:saddles}(a), which is single-valued around $t_i$. Hence
$c_i(\Delta)$ is this value. Equation~\eqref{eq:affine} is affine in
$(W,\Delta)$, and its homogeneous part does not depend on $\Delta$; by the
uniqueness in Lemma~\ref{lem:saddles}(a), $\gamma_i$ depends linearly on
$\Delta$, and so does $c_i(\Delta)$. The linear form
$\Delta\mapsto c_2(\Delta)-c_1(\Delta)$ vanishes exactly when $f_1$ and
$f_2$ have a common fixed point, that is, by Lemma~\ref{lem:polsol}, on
the kernel of $\lambda$. Hence $c_2-c_1=\kappa_0\lambda$ with
$\kappa_0\neq0$. Put $\kappa=(1-\zeta)\kappa_0$. In the coordinate $W'$,
the map $f_1$ is the rotation of order three about $0$, that is,
$W'\mapsto\zeta W'$, and $f_2$ is the rotation of order three about
$\kappa_0\lambda(\Delta)/\kappa=\lambda(\Delta)/(1-\zeta)$, that is,
$W'\mapsto\zeta W'+\lambda(\Delta)$. In the coordinate $x'=1/W'$, these
maps are $x'\mapsto\zeta^2x'$ and $x'\mapsto x'/(\zeta+\lambda(\Delta)
x')$, and the change of coordinate $x\mapsto x'$, where $x=1/W$, is the
germ of biholomorphism $x'=\kappa x/(1-c_1(\Delta)x)$.
\end{proof}

\begin{proof}[Proof of Theorem~\textup{\ref{thm:group}}]
We use the coordinate $W=1/x$, in which $g_1(W)=\zeta W$ and
$g_{2,\lambda}(W)=\zeta W+\lambda$. These are Möbius transformations of the
Riemann sphere fixing $\infty$, that is, $x=0$.

(a) Both maps are affine with linear part $\zeta$, and
$g_{2,\lambda}^3(W)=W+\lambda(1+\zeta+\zeta^2)=W$.

(d) The map $\tau=g_{2,\lambda}\circ g_1^{-1}$ is the translation
$W\mapsto W+\lambda$, and $g_1^j\circ\tau\circ g_1^{-j}$ is the
translation by $\zeta^j\lambda$. Let $\lambda\neq0$ and
$\Lambda=\Z\lambda+\Z\zeta\lambda$, a lattice since $\zeta\notin\R$, with
$\zeta\Lambda=\Lambda$ because $\zeta^2=-1-\zeta$. The maps
$W\mapsto\zeta^jW+\mu$, $j\in\Z/3\Z$, $\mu\in\Lambda$, form a group, which
contains $g_1$ and $g_{2,\lambda}=\tau\circ g_1$ and is generated by $g_1$
and the conjugates of $\tau$. Hence it is $G_\lambda$, which is therefore
isomorphic to $\Z^2\rtimes\Z/3\Z$. The homothety $H_\lambda(x)=\lambda x$
satisfies $H_\lambda\circ g_1\circ H_\lambda^{-1}=g_1$ and
$H_\lambda\circ g_{2,\lambda}\circ H_\lambda^{-1}=g_{2,1}$, so that
$H_\lambda G_\lambda H_\lambda^{-1}=G_1$. Finally,
$g_1\circ g_{2,\lambda}^{-1}$ is the translation $W\mapsto W-\lambda$,
that is, $x\mapsto x/(1-\lambda x)$, which is tangent to the identity and
different from it if $\lambda\neq0$.

(b) Each element of the pseudogroup on $D_r$ is the restriction of an
element of $G_\lambda$, so the orbit of a point is contained in its
$G_\lambda$-orbit. In the coordinate $W$, the $G_\lambda$-orbit of $W_0$
is contained in the finite set $\{\zeta^jW_0\}$ if $\lambda=0$, and in the
discrete set $\bigcup_j(\zeta^jW_0+\Lambda)$ if $\lambda\neq0$. In both
cases it meets each annulus $\{\epsilon\leq|x|\leq r'\}$, $0<\epsilon<r'<r$,
in finitely many points, so that every orbit of the pseudogroup through
a non-zero point is a closed subset of $D_r\setminus\{0\}$. If $\lambda=0$, the function $x^3$ is
invariant. If $\lambda\neq0$, the elliptic function
$\mathcal R(W)=\sum_{\mu\in\Lambda}(W-\mu)^{-3}$, a multiple of the
derivative of the Weierstrass function of $\Lambda$, is invariant under
the translations by $\Lambda$ and, since $\zeta\Lambda=\Lambda$ and
$\zeta^3=1$, under $W\mapsto\zeta W$; hence $x\mapsto\mathcal R(1/x)$ is
a non-constant $G_\lambda$-invariant meromorphic function on
$\C\setminus\{0\}$.

(c) If $\lambda=0$, then $g_{2,0}=g_1$, so $G_0$ is the cyclic group of
order three generated by $x\mapsto\zeta^2x$, and $x^3$ is invariant. If
$\lambda\neq0$, the parabolic germ $g_1\circ g_{2,\lambda}^{-1}$ has
infinite order, so $G_\lambda$ is infinite, and a non-constant
$G_\lambda$-invariant germ would have an infinite invariance group,
contradicting Lemma~\ref{lem:finite}(b).

The last assertion of Theorem~\ref{thm:group} is Lemma~\ref{lem:normal}.
\end{proof}

\subsection{The multivalued first integral and the closed leaves}
\label{sec:Phi}

Let $T=\C\setminus\{0,1\}$, let $\pi_T\colon\widetilde T\to T$ be the
universal covering, and fix on $\widetilde T$ a holomorphic cube root
$\varrho$ of $P\circ\pi_T$ and a primitive $J$ of the holomorphic $1$-form
$\varrho^{-1}\pi_T^*\bigl(\delta(t)\,dt/P(t)\bigr)$. Let $U$ be a connected
open neighbourhood of the origin in which $\omega_\Delta$ has no zero
other than $0$, let $\Fcal_U$ be the corresponding representative of
$\Fcal_\Delta$, and $U^*=U\setminus\{F=0\}$. The set
$\{F=0\}\cap U\setminus\{0\}$ is a union of leaves of $\Fcal_U$, so every
other leaf is contained in $U^*$. On $U^*$ we have $x\neq0$ and
$t=y/x\in T$. Let $\widetilde U=U^*\times_T\widetilde T$ be the fibred
product, and $\pi_U\colon\widetilde U\to U^*$ the projection, which is a
covering map; the group $\mathrm{Deck}$ of deck transformations of $\pi_T$
acts on $\widetilde U$ through the second factor, transitively on the
fibres of $\pi_U$. We consider the function
\[
\Phi=\frac1{x\,\varrho}-J\colon\widetilde U\to\C ;
\]
in other words, $\Phi$ is the multivalued function $F^{-1/3}-J(y/x)$ on
$U^*$.

\begin{lemma}\label{lem:Phi}
\begin{enumerate}[label=\textup{(\alph*)}]
\item On $\widetilde U$ we have $(x\varrho)^3=F$ and
\[
d\Phi=-\frac{\omega_\Delta}{x\,F\,\varrho},
\]
where functions and forms on $U^*$ are identified with their pull-backs
by $\pi_U$. In particular, $\Phi$ is constant along the lift to
$\widetilde U$ of any path contained in a leaf of\/ $\Fcal_U$, and its
restriction to any lift of a transverse section contained in $U^*$ is a
local biholomorphism.
\item For every $\gamma\in\mathrm{Deck}$ there is an affine map
$g_\gamma(u)=\chi(\gamma)^{-1}u-c(\gamma)$, with $\chi(\gamma)^3=1$, such
that $\Phi\circ\gamma=g_\gamma\circ\Phi$. The maps $g_\gamma$ form a group
$G_\Phi$, which is generated by two elements $B_1$ and $B_2$ with linear
part $\zeta^{-1}$.
\item There is a non-constant meromorphic function $\mathcal R$ on $\C$
such that $\mathcal R\circ g=\mathcal R$ for every $g\in G_\Phi$.
\end{enumerate}
\end{lemma}

\begin{proof}
(a) The first identity follows from $x^3P(y/x)=F$. Since
$\varrho^{-1}d\varrho=\frac13\,(P'/P)\,dt$ and
$dJ=\delta\,dt/(P\varrho)$, we get, in the coordinates $(t,W)$,
\[
d\Phi=\frac1\varrho\Bigl(dW-\frac{P'}{3P}\,W\,dt-\frac{\delta}{P}\,dt\Bigr)
=-\frac{1}{x^2P\varrho}\Bigl(P\,dx+x\bigl(\tfrac13\,P'+x\,\delta\bigr)\,dt\Bigr)
=-\frac{\pi^*\omega_\Delta}{x^4P\varrho},
\]
and $x^4P(y/x)=xF$. The tangent vectors of the lift of a path contained in
a leaf are annihilated by $\omega_\Delta$, hence by $d\Phi$, and
$\omega_\Delta$ vanishes nowhere on a transverse section.

(b) The function $\varrho\circ\gamma$ is another cube root of
$P\circ\pi_T$, so $\varrho\circ\gamma=\chi(\gamma)\varrho$ with
$\chi(\gamma)^3=1$, and then $d(J\circ\gamma)=\chi(\gamma)^{-1}dJ$, so that
$J\circ\gamma=\chi(\gamma)^{-1}J+c(\gamma)$ for a constant $c(\gamma)$;
hence $\Phi\circ\gamma=\chi(\gamma)^{-1}\Phi-c(\gamma)$. Since $\Phi$ is an
open map, $\Phi\circ\gamma\gamma'=g_\gamma\circ g_{\gamma'}\circ\Phi$
implies $g_{\gamma\gamma'}=g_\gamma g_{\gamma'}$. The group
$\mathrm{Deck}$ is free on two generators which correspond to the positive
loops around $0$ and $1$, along which $\varrho$ is multiplied by $\zeta$.

(c) By (b), $B_1$ and $B_2$ are rotations of $\C$, of order three, about
two points $c'_1$ and $c'_2$. If $c'_1=c'_2$, then $G_\Phi$ is generated
by $B_1$, and we put $\mathcal R(u)=(u-c'_1)^3$. Otherwise, in the
coordinate $v=(u-c'_1)/\kappa$, $\kappa=(1-\zeta^{-1})(c'_2-c'_1)$, we
have $B_1(v)=\zeta^2v$ and $B_2(v)=\zeta^2v+1$, so that $G_\Phi$ becomes
the group of the maps $\zeta^jv+\mu$, $j\in\Z/3\Z$,
$\mu\in\Lambda_1=\Z+\Z\zeta$, which is the group $G_1$ of
Theorem~\ref{thm:group} in its coordinate $W$ (see the proof of its
part~(d)). We put $\mathcal R(u)=\sum_{\mu\in\Lambda_1}(v-\mu)^{-3}$, the
$G_1$-invariant function of the proof of Theorem~\ref{thm:group}(b). It
is an elliptic function of $u$ with respect to the lattice
$\Lambda_\Phi=\kappa\Lambda_1$, with a single pole, of order three, in
each period parallelogram.
\end{proof}

\begin{proof}[Proof of Theorem~\textup{\ref{thm:Q52}(b)}]
By Lemma~\ref{lem:Phi}(b),(c), the map $\mathcal R\circ\Phi$ is invariant
under $\mathrm{Deck}$, so it descends to a holomorphic map
$\mathcal R_\Phi\colon U^*\to\Pp^1$. By Lemma~\ref{lem:Phi}(a),
$\mathcal R_\Phi$ is constant on the leaves contained in $U^*$, and its
restriction to a transverse section contained in $U^*$ is locally the
composition of $\mathcal R$ with a biholomorphism onto an open subset of
$\C$; hence this restriction is non-constant.

Let $L\subset U^*$ be a leaf and $z\in U^*$ a point of its closure.
Choose a flow box $V\subset U^*$ at $z$ and a transverse section
$\Sigma\subset V$ at $z$ which meets every plaque of $V$ and on which
$\mathcal R_\Phi$ takes the value $\mathcal R_\Phi(z)$ only at $z$. By
continuity $\mathcal R_\Phi\equiv\mathcal R_\Phi(z)$ on $L$, so every
plaque of $V$ contained in $L$ meets $\Sigma$ at a point where
$\mathcal R_\Phi=\mathcal R_\Phi(z)$, that is, at $z$. Hence $L\cap V$ is
the plaque of $z$, and $z\in L$. Therefore $L$ is closed in $U^*$.
\end{proof}

Thus, when $\lambda(\Delta)\neq0$, the germ $\Fcal_\Delta$ is a
non-dicritical generalized curve whose leaves are closed off the
separatrices and whose three separatrices have trivial, hence periodic,
holonomy, but which has no non-constant holomorphic first integral.

By Lemma~\ref{lem:Phi}(a), $\Phi$ is a primitive of the closed form
$-\omega_\Delta/F^{4/3}$, so that it is a Liouvillian first integral of
$\Fcal_\Delta$, in accordance with the third alternative of
\cite[Theorem~1.2]{Sca16}.

\subsection{Growth of the leaf counts}\label{sec:growthC}

The exponent in Theorem~\ref{thm:Q52}(e) has a simple origin. Transverse
to the exceptional line, the trace of a leaf is essentially an orbit of
the lattice $\Lambda$ of translations of the projective holonomy group, so
that the counts grow like $\rho^{-2}$, the exponent $2$ being the rank of
$\Lambda$; the passage through the corner $s_i$, where $\tau\xi^3$ is a
first integral, divides this exponent by $3$.

\begin{proof}[Proof of Theorem~\textup{\ref{thm:Q52}(e)}]
Let $\lambda(\Delta)\neq0$, let $\Gamma=\Gamma_i$, let $\Fcal_V$ be a
representative adapted to $\Gamma$, and let $S$ be a sector transverse to
$\Gamma$, of opening $\vartheta$, in a transverse section $\Sigma$ at a
point $a\in\Gamma^\circ_V$, with coordinate $z$. We use the blow-up $\pi$
over $V$; for transverse sections contained in $\pi^{-1}(V)\setminus E$,
the leaves of $\tF_\Delta$ on $\pi^{-1}(V)$ correspond to those of
$\Fcal_V$, and we identify such sections with their images. The function
$\Phi$ of Section~\ref{sec:Phi} and the function $\mathcal R_\Phi$ of the
proof of Theorem~\ref{thm:Q52}(b) are defined, with the same proofs, on
$V\setminus\{F=0\}$: the hypothesis on the zeros of $\omega_\Delta$ is not
used there, and $\mathcal R_\Phi$ is constant on the leaves of $\Fcal_V$
contained in $V\setminus\{F=0\}$.

\emph{Step~1: counts transverse to $E$.} Let $\Sigma_*$ be the transverse
section $\{t=t_*\}$ to $E$ at the point $(t_*,0)$ of the chart $y=tx$,
with the coordinate $x$, and $D_*=\{0<|x|<r_*\}\subset\Sigma_*$, with
$r_*$ so small that $D_*\subset\pi^{-1}(V)$. We show:
\begin{enumerate}[label=\textup{(\roman*)}]
\item there is $C_*$ such that every leaf meets $D_*\cap\{|x|>\rho\}$ in
at most $C_*\rho^{-2}$ points, for $0<\rho<1$;
\item for every sector $S_*$ of $\Sigma_*$ there are $c_*,\rho_*>0$ with
$\nu_V(S_*,\rho)\geq c_*\rho^{-2}$ for $0<\rho<\rho_*$.
\end{enumerate}

(i) All the points of $D_*$ have the same coordinate $t=t_*$, so $D_*$
lifts to the fibred product $\widetilde U$ of Section~\ref{sec:Phi} by
fixing a point $\tilde t_*$ over $t_*$; on this lift
$\Phi=\alpha W+\beta$ with $W=1/x$, $\alpha=1/\varrho(\tilde t_*)\neq0$
and $\beta=-J(\tilde t_*)$. Hence the trace of a leaf on $D_*$ is
contained in a set $\{W:\ \mathcal R(\alpha W+\beta)=c\}$, $c\in\Pp^1$.
The function $\mathcal R$ of Lemma~\ref{lem:Phi}(c) is elliptic with
respect to the lattice $\Lambda_\Phi$ of its proof, with a single pole, of
order three, in each period parallelogram; it therefore takes each value
exactly three times, counted with multiplicity, in each period
parallelogram, and $\{\mathcal R=c\}$ is contained in the union of at
most three cosets $u_k+\Lambda_\Phi$. A coset of $\Lambda_\Phi$ has at
most $C(1+A)^2$ points in a disc of radius $A$. The points with
$|x|>\rho$ have $|\alpha W+\beta|<|\alpha|/\rho+|\beta|$, which
gives~(i).

(ii) We use the coordinate $W'=(W-c_1)/\kappa$ of
Lemma~\ref{lem:normal}, with $c_1=c_1(\Delta)$ and
$\kappa=(1-\zeta)\kappa_0$, in which the holonomy maps $f_1$ and $f_2$
are $\zeta W'$ and $\zeta W'+\lambda$, $\lambda=\lambda(\Delta)$. Let
$\ell_0$ and $\ell_1$ be loops in $T=\C\setminus\{0,1\}$ based at $t_*$
whose holonomy maps are $f_2\circ f_1^{-1}$ and
$f_1\circ(f_2\circ f_1^{-1})\circ f_1^{-1}$, the translations by
$\lambda$ and by $\zeta\lambda$. Over a path $c$ in $T$, the lifts of $c$
to the leaves of $\tF_\Delta$ off $E$ are the graphs of the solutions
of~\eqref{eq:affine} along $c$, and continuation along $c$ is an affine
map $W\mapsto A_c(t)W+B_c(t)$ with $A_c(t)\neq0$ depending continuously on
$t\in c$. Since $\pi^{-1}(V)$ is a neighbourhood of $E$, it contains
$\{(t,x):\ t\in K,\ |x|<r_K\}$ for every compact set $K\subset\C$ of the
chart. Hence there is $R_1>0$ such that, for every point of $\Sigma_*$
with $|W'|\geq R_1$, the lifts of $\ell_0^{\pm1}$ and $\ell_1^{\pm1}$
starting at this point stay in $\pi^{-1}(V)\setminus E$; their end points,
$W'\pm\lambda$ and $W'\pm\zeta\lambda$, lie on the same leaf of
$\tF_\Delta$ on $\pi^{-1}(V)$. Let $\mu=m\lambda+n\zeta\lambda\in\Lambda$,
and go from $W'_0$ to $W'_0+\mu$ by $|m|$ steps $\pm\lambda$ followed by
$|n|$ steps $\pm\zeta\lambda$. The intermediate points are
$W'_0+\nu$ with $|\nu|\leq(|m|+|n|)|\lambda|\leq C_0|\mu|$, where $C_0$
depends only on the basis $(\lambda,\zeta\lambda)$. Hence
\emph{$W'_0+\mu$ lies on the leaf of $W'_0$ as soon as
$|W'_0|-C_0|\mu|\geq R_1$.}

Let $S_*=S(r;\theta_1,\theta_1+\vartheta_*)$ in the coordinate $x$,
$\theta_c=\theta_1+\vartheta_*/2$, and
$\epsilon=\min\{\sin(\vartheta_*/4),1/(4C_0)\}$. For $R>0$ let $W'_0$ be
defined by $\kappa W'_0+c_1=|\kappa|R\,e^{-i\theta_c}$. For $\mu\in\Lambda$
with $|\mu|\leq\epsilon R$, the point $x=1/(\kappa(W'_0+\mu)+c_1)$
satisfies
\[
\frac1{2|\kappa|R}<\frac1{(1+\epsilon)|\kappa|R}\leq|x|
\leq\frac1{(1-\epsilon)|\kappa|R},\qquad|\arg x-\theta_c|\leq
\arcsin\epsilon\leq\frac{\vartheta_*}4 ,
\]
because $1/x=|\kappa|R\,e^{-i\theta_c}+\kappa\mu$ and
$|\kappa\mu|\leq\epsilon|\kappa|R$; so $x\in S_*$ for $R$ large; moreover
$|W'_0|-C_0|\mu|\geq R-|c_1|/|\kappa|-R/4\geq R_1$ for $R$ large. These
points are pairwise distinct, lie on one leaf, and their number is at
least $c_\Lambda\epsilon^2R^2$ for $R$ large. With $\rho=1/(2|\kappa|R)$
this gives~(ii).

\emph{Step~2: transport to $\Gamma$.} Let $(\xi,\tau)$ be the coordinates
of Lemma~\ref{lem:saddles}(b) at $s_i$, and $W_i\subset\pi^{-1}(V)$ a
neighbourhood of $s_i$ as in Lemma~\ref{lem:connected} for the first
integral $g=\tau\xi^3$, which is not a power. Choose $\xi_0\neq0$ and
$\tau_0\neq0$ so small that the transverse sections
$\Sigma_1=\{\xi=\xi_0\}$ to the strict transform of $\Gamma$, with the
coordinate $\tau$, and $\Sigma_E=\{\tau=\tau_0\}$ to $E$, with the
coordinate $\xi$, contain discs about $(\xi_0,0)$ and $(0,\tau_0)$
contained in $W_i$. Points $(\xi,\tau_0)$ and $(\xi_0,\tau)$ of these
discs with $\tau_0\xi^3=\xi_0^3\tau\neq0$ have the same non-zero value of
$g$, hence lie on a common leaf, by Lemma~\ref{lem:connected}(a). Let
$\sigma_\Gamma\colon(\Sigma,a)\to(\Sigma_1,(\xi_0,0))$ and
$\sigma_E\colon(\Sigma_E,(0,\tau_0))\to(\Sigma_*,(t_*,0))$ be the
holonomy maps of paths in the leaves $\Gamma^\circ_V$ and $E^\circ$; they
are germs of biholomorphisms, so that
$c|z|\leq|\tau(\sigma_\Gamma(z))|\leq C|z|$ and
$c|\xi|\leq|x(\sigma_E(\xi))|\leq C|\xi|$ near the base points, for
constants $c,C>0$. We use repeatedly the injection argument of
Lemma~\ref{lem:counts}(c): an injective map between subsets of transverse
sections which preserves the leaves and satisfies $|w'|\geq c|w|^\kappa$
bounds the number of points of a leaf in the source, outside radius
$\rho$, by the number in the target outside radius $c\rho^\kappa$.

\emph{Lower bound.} By Lemma~\ref{lem:counts}(d), applied to
$\sigma_\Gamma$, there are a sector $S_1$ of $\Sigma_1$, in the coordinate
$\tau$, and $c_1>0$ with $\nu_V(S_1,\rho)\leq\nu_V(S,c_1\rho)$. Choose an
open interval $I$ of length less than $2\pi/3$ such that
$\arg(\tau_0/\xi_0^3)+3I$ is compactly contained in the interval of
arguments of $S_1$. On a sector $S_E$ of $\Sigma_E$ of arguments $I$ and
small radius, the map $\beta(\xi)=\tau_0\xi^3/\xi_0^3$ is injective, takes
its values in $S_1$ and satisfies $|\beta(\xi)|=c_2|\xi|^3$; hence
$\nu_V(S_E,\rho)\leq\nu_V(S_1,c_2\rho^3)$. By Lemma~\ref{lem:counts}(d),
applied to $\sigma_E^{-1}$, there are a sector $S_*$ of $\Sigma_*$, in the
coordinate $x$, and $c_3>0$ with $\nu_V(S_*,\rho)\leq\nu_V(S_E,c_3\rho)$.
With Step~1(ii), for small $\rho$,
\[
c_*\rho^{-2}\leq\nu_V(S_*,\rho)\leq\nu_V(S_E,c_3\rho)
\leq\nu_V(S_1,c_2c_3^3\rho^3)\leq\nu_V(S,c_1c_2c_3^3\rho^3),
\]
that is, $\nu_V(S,\rho')\geq c\,\rho'^{-2/3}$ for small $\rho'$.

\emph{Upper bound.} If the radius of $S$ is small, $\sigma_\Gamma$ maps
$S$ injectively into a punctured disc $D_1=\{0<|\tau|<r_1\}$ of
$\Sigma_1$, with $|\tau(\sigma_\Gamma(z))|\geq c|z|$, and $r_1$ is as
small as we wish. The punctured disc $D_1$ is the union of the two sectors
$D_1^\pm$ of opening $2\pi$ whose arguments lie in $(0,2\pi)$ and in
$(-\pi,\pi)$. On each of them a branch of
$\beta^\pm(\tau)=\xi_0(\tau/\tau_0)^{1/3}$ is injective, preserves the
leaves, and satisfies $|\beta^\pm(\tau)|=c_4|\tau|^{1/3}$; for $r_1$
small its image lies in a punctured disc of $\Sigma_E$ on which
$\sigma_E$ is defined and injective, with values in $D_*$. Hence, for a
leaf $L$,
\[
\#\bigl(L\cap S\cap\{|z|>\rho\}\bigr)
\leq\sum_\pm\#\bigl(L\cap D_1^\pm\cap\{|\tau|>c\rho\}\bigr)
\leq2\,\#\bigl(L\cap D_*\cap\{|x|>cc_4(c\rho)^{1/3}\}\bigr),
\]
and Step~1(i) gives $\nu_V(S,\rho)\leq C\rho^{-2/3}$ for small $\rho$.
\end{proof}

\section{Analytic graphs: proof of Theorem~\ref{thm:graph}}\label{sec:graph}

Let $\Fcal_U$ be a representative of a germ $\Fcal$ at $0\in\C^2$, and
$\mathcal G=\mathcal G(\Fcal_U)$ its graph (Definition~\ref{def:graph}).

\begin{lemma}\label{lem:saturated}
If $(p,q)\in\mathcal G$, with $p,q\notin\sing(\Fcal)$, and $p'$, $q'$ lie
on the leaves of $\Fcal_U$ through $p$ and $q$, then
$(p',q')\in\mathcal G$. Consequently, if $\mathcal G\neq U\times U$ is
analytic, then for every transverse section $\Sigma$ the set
$\mathcal G\cap(\Sigma\times\Sigma)$ has empty interior in
$\Sigma\times\Sigma$.
\end{lemma}

\begin{proof}
Let $c$ be a path in the leaf of $p$ from $p$ to $p'$, covered by a chain
of flow boxes. Let $(p_k,q_k)\to(p,q)$ with $p_k$ and $q_k$ on a common
leaf. For $k$ large, following the plaques of the chain from $p_k$ gives a
point $p'_k$ on the leaf of $p_k$, hence of $q_k$, with $p'_k\to p'$.
Thus $(p'_k,q_k)\to(p',q)$, and $(p',q)\in\mathcal G$; the same argument
in the second factor gives $(p',q')\in\mathcal G$. If
$\mathcal G\cap(\Sigma\times\Sigma)$ contained a product $W_1\times W_2$ of
nonempty open subsets of $\Sigma$, the first assertion would give
$\mathcal G\supset\mathrm{Sat}(W_1)\times\mathrm{Sat}(W_2)$, where
$\mathrm{Sat}(W_i)$ is the union of the leaves meeting $W_i$, an open
subset of $U$. An analytic subset of the connected manifold $U\times U$
with nonempty interior is $U\times U$.
\end{proof}

\begin{proposition}\label{prop:graph-counts}
Suppose that $\mathcal G=\mathcal G(\Fcal_U)\neq U\times U$ is analytic.
Let $a\in U\setminus\sing(\Fcal)$ and let $\Sigma$ be a transverse section
at $a$, with coordinate $z$. Then there are $r>0$ and $K\geq1$ such that
every leaf of $\Fcal_U$ meets the punctured disc $\{0<|z|<r\}$ in at most
$K$ points.
\end{proposition}

\begin{proof}
The set $Z=\mathcal G\cap(\Sigma\times\Sigma)$ is an analytic subset of
$\Sigma\times\Sigma$, which contains the pairs of points of $\Sigma$ on a
common leaf; by Lemma~\ref{lem:saturated} it has dimension at most one.
Let $\overline D\subset\Sigma$ be a closed disc about $a$. Only finitely
many irreducible components of $Z$ meet the compact set
$\overline D\times\overline D$. Call a component \emph{vertical} if it is
contained in $\{c\}\times\Sigma$ for some $c$; there are finitely many
such points $c$, and we choose $r$ so small that none of them lies in
$\{0<|z|<r\}$. Each of the other components $Z_j$ meets
$\{a\}\times\overline D$ in finitely many points $(a,w_{j,l})$. Near each
of them, the projection $Z_j\to\Sigma$ onto the first factor is a finite
branched covering over a neighbourhood of $a$, of some degree $d_{j,l}$.
Since $Z_j\cap(\Sigma\times\overline D)$ is closed, for $|z|$ small all
the points of $Z_j$ over $z$ in $\Sigma\times\overline D$ lie near the
points $(a,w_{j,l})$. Hence, decreasing $r$, every $z$ with $0<|z|<r$ has
at most $K=\sum_{j,l}d_{j,l}$ points $w\in\overline D$ with $(z,w)\in Z$.
If a leaf meets $\{0<|z|<r\}$ at $z$, all its points in this punctured
disc are such points $w$.
\end{proof}

\begin{proof}[Proof of Theorem~\textup{\ref{thm:graph}}]
Suppose first that $\Fcal$ has a non-constant holomorphic first integral.
By Lemma~\ref{lem:integrable-gc}, $\Fcal$ is a non-dicritical
generalized curve. Let $F_0$ be a primitive first integral; by
Lemma~\ref{lem:primitive}, $g=F_0-F_0(0)$ is not a power. In a given
neighbourhood of the origin, let $V$ be a connected open neighbourhood of
$0$ as in Lemma~\ref{lem:connected} for $g$, on which $0$ is the only
singular point. Let $p,q\in V\setminus\{0\}$. If $p$ and $q$ lie on a
common leaf of $\Fcal_V$, then $F_0(p)=F_0(q)$. Conversely, if
$F_0(p)=F_0(q)\neq F_0(0)$, then $p$ and $q$ lie on a common leaf, by
Lemma~\ref{lem:connected}(a). Thus the set of pairs on a common leaf is
contained in the hypersurface $A=\{F_0(p)=F_0(q)\}$ of $V\times V$, and
contains $A\cap\{F_0(p)\neq F_0(0)\}$. The complement of the latter set
in $A$ is $\{F_0(p)=F_0(q)=F_0(0)\}$, which has dimension two and
therefore contains no irreducible component of $A$, all of which have
dimension three. Hence $\mathcal G(\Fcal_V)=A$, an analytic subset of
$V\times V$ different from $V\times V$.

Conversely, let $\Fcal$ be non-dicritical, and let $\Fcal_U$ be a
representative whose graph is analytic and different from $U\times U$.
Let $\Gamma$ be a separatrix of $\Fcal$. Choose a representative
$\Fcal_{U'}$ adapted to $\Gamma$, with $U'\subset U$
(Lemma~\ref{lem:sep}), a point $a\in\Gamma^\circ_{U'}$, and a transverse
section $\Sigma\subset U'$ at $a$. By
Proposition~\ref{prop:graph-counts}, small sectors of $\Sigma$ with vertex
$a$ are met by the leaves of $\Fcal_U$ in a uniformly bounded number of
points, hence also by the leaves of $\Fcal_{U'}$
(Lemma~\ref{lem:counts}(b)). Bounded counts are sub-logarithmic, and this
holds for every separatrix, so Theorem~\ref{thm:plane}(a) gives a
non-constant holomorphic first integral.
\end{proof}

The hypothesis that $\Fcal$ be non-dicritical cannot be dropped: the radial
foliation $x\,dy-y\,dx$ has the analytic graph $\{x_1y_2-x_2y_1=0\}$ and no
non-constant holomorphic first integral. In one variable, bounded counts
do not make the graph of a relation analytic: identifying each point $z$
of the disc of centre $4^{-k}$ and radius $4^{-k}/10$ with $z-4^{-k}/2$,
$k\geq1$, gives an open relation with classes of at most two points whose
graph contains pieces of the infinitely many lines $w=z-4^{-k}/2$
accumulating at $(0,0)$. In Theorem~\ref{thm:graph}, the passage from
bounded counts to an analytic graph is provided by the first integral of
Theorem~\ref{thm:plane}.

\bigskip
\noindent{\sffamily\bfseries AI disclosure}\enspace
The mathematical ideas, results and proofs presented in this paper were
conceived and developed by the authors. Claude (Anthropic) was used as an
auxiliary tool to assist in checking calculations, examining mathematical
arguments, and drafting and revising the exposition. The authors take full
responsibility for the mathematical correctness and content of the paper.

\end{document}